%% file: main.tex
\documentclass[letterpaper]{article} 
\usepackage{aaai24}  
\usepackage{times}  
\usepackage{helvet}  
\usepackage{courier}  
\usepackage[hyphens]{url}  
\usepackage{graphicx} 
\usepackage{natbib}  
\usepackage{caption} 
\usepackage{algorithm}
\usepackage{algorithmic}

\usepackage{newfloat}
\usepackage{listings}
\DeclareCaptionStyle{ruled}{labelfont=normalfont,labelsep=colon,strut=off} 
\floatstyle{ruled}
\newfloat{listing}{tb}{lst}{}
\floatname{listing}{Listing}
\title{Decision-Making for Land Conservation: A Derivative-Free Optimization Framework with Nonlinear Inputs\footnote{This is the extended version of our paper that was accepted to AAAI'24. This version includes the technical appendix.}}
\author{
    Cassidy K. Buhler, Hande Y. Benson
}
\affiliations{
    Department of Decision Sciences and MIS, Drexel University\\

    cb3452@drexel.edu, 
    hvb22@drexel.edu\\
}
\nocopyright

\usepackage{bibentry}

\usepackage{amsmath}
\usepackage{amssymb}
\usepackage[T1]{fontenc}
\usepackage{xcolor}

\usepackage{appendix}
\begin{document}

\maketitle

\begin{abstract}

Protected areas (PAs) are designated spaces where human activities are restricted to preserve critical habitats. Decision-makers are challenged with balancing a trade-off of financial feasibility with ecological benefit when establishing PAs. Given the long-term ramifications of these decisions and the constantly shifting environment, it is crucial that PAs are carefully selected with long-term viability in mind. 

Using AI tools like simulation and optimization is common for designating PAs, but current decision models are primarily linear. In this paper, we propose a derivative-free optimization framework paired with a nonlinear component, population viability analysis (PVA). Formulated as a mixed integer nonlinear programming (MINLP) problem, our model allows for linear and nonlinear inputs. Connectivity, competition, crowding, and other similar concerns are handled by the PVA software, rather than expressed as constraints of the optimization model. In addition, we present numerical results that serve as a proof of concept, showing our models yield PAs with similar expected risk to that of preserving every parcel in a habitat, but at a significantly lower cost. 

The overall goal is to promote interdisciplinary work by providing a new mathematical programming tool for conservationists that allows for nonlinear inputs and can be paired with existing ecological software.
Our code and data are available at
https://github.com/cassiebuhler/conservation-dfo.

\end{abstract}

\section{Introduction}

As a consequence of habitat degradation due to human activities, protected areas (PAs) have been implemented for protecting vulnerable species and preserving valuable ecological processes  \cite{world1980world}. In recent years, the post-2020 Global Biodiversity Framework has mandated that 30\% of land and sea be established as protected areas by 2030 \cite{GDF2020,EO_30}. This initiative, coined ``$30\times30$'', has motivated further research on PAs and their decision frameworks. 

Conservationists have cautioned that area-based goals, such as $30\times30$, have the potential to prioritize quantity over quality and emphasize the importance of prioritizing ecological value, such as biodiversity, in these decisions \cite{maxwell2020area, joppa2009high}. However, financial support is still vital to the long-term management and efficacy of a PA 
\cite{coad2019widespread,rodrigues2020multifaceted, watson2014performance, lawler2020planning}.

When implemented appropriately, PAs are a necessary investment for our planet \cite{dudley2008guidelines}. With the environment rapidly changing (e.g. climate change increasing habitat fragmentation \cite{costanza2020preserving}), it is crucial that PAs are carefully selected with consideration of future conditions to reduce the risk of extinction \cite{williams2022global,akasaka2017smart,cruz2021global}.

\subsection{Population Viability Analysis}
While there are many ways to estimate extinction risk, one of the more popular tools in conservation biology is population viability analysis (PVA), which predicts the probability of extinction for a species in a particular habitat  \cite{akccakaya2000population}. Scientists use forecasts from PVA in recovery plans  \cite{runge2017status, faust2016red} and management decisions \cite{lacy2021population} for threatened species. 

PVA also provides insight on the persistence of a species in PAs, ultimately used to inform practitioners on improvements necessary to increase the likelihood of survival \cite{winton2020protected, finnegan2021reserve}. In addition, PVA has been used to advise potential PA locations based on an unprotected population that is declining \cite{andersen2023influence}. 

Researchers have noted that PVA is a valuable risk assessment, yet the necessary data required for PVA can be challenging to obtain, which often limits PVA to single-species \cite{rondinini2010quantitative}. Moreover, many PVA studies in the literature were found to lack clear documentation of parameters and assumptions \cite{doak2015recommendations}
making it impossible to reproduce results \cite{morrison2016repeatability} which ultimately compromises its quality and reliability \cite{chaudhary2020critical} as code and data transparency is necessary in conservation decision-making  \cite{morrison2016repeatability}.

\subsection{Mathematical Programming}
Establishing PAs requires balancing a trade-off of financial feasibility with ecological benefit.
For decades, optimization models have aided decision-makers with such trade-offs \cite{alagador2022operations, billionnet2013mathematical}. 
The first use of mathematical programming to select conservation sites was published by \cite{margules1988selecting} and \cite{cocks1989using}. The authors proposed using integer programming (IP), with linear inputs and binary decisions, rather than existing ranking methods. 

In mathematical programming models. a landscape is modeled as a raster divided into $n\times n$ pixels, where each pixel is a parcel eligible to be selected as a PA. 
The set of parcels is denoted as $P$, and the set of species is denoted as $S$. Let $x_p, p \in P$ be the binary decision variable denoting whether or not parcel $p$ on a landscape is selected as a PA. The IP is defined as follows:
 \begin{align} \label{minIP}
    \displaystyle \underset{x  \in  \{0,1\}^{|P|} }{\textrm{minimize}} \quad c^T x \quad \textrm{subject to} \quad Ax \geq  b,
\end{align}
where $A \in \mathbb{R}^{|S|\times |P|}$, $c \in \mathbb{R}^{|P|}$, and $b \in \mathbb{R}^{|S|}$.  In the minimum set cover formulation, $c_i$ is the cost of acquisition for parcel $i$, $A_{j,i}$ is the population for species $j$ in parcel $i$, and $b_j$ is a target population for species $j$. 
Following the introduction of optimization models in conservation, variations of \eqref{minIP} were proposed 
\cite{underhill1994optimal,possingham1993mathematics,revelle2002counterpart}.  IPs also can represent the maximal set coverage: maximizing a species population without exceeding a budget $\beta$. 
 \begin{align}
    \displaystyle \underset{x  \in  \{0,1\}^{|P|} }{\textrm{maximize}} \quad \sum_{i \in P, j \in S} A_{ji} x_i  \quad \textrm{subject to} \quad  c^T x \leq \beta
\end{align}
Maximal set coverage is also used in many applications \cite{camm1996note,church1996reserve,arthur1997finding,csuti1997comparison,rosing2002maximizing,rodrigues2000flexibility}.  More recently, spatial properties, such as compactness  \cite{marianov2008selecting} and functional connectivity \cite{onal2016optimal,dilkina2017trade,gupta2019reserve, costanza2020preserving, williams2020incorporating}, have been introduced as constraints in these models. 

To this day, IP is the most common optimization  model used in conservation planning. This formulation limits the objective and constraints to linearity, which ultimately, limits the potential applications of the model.

\subsubsection{Mixed Integer Linear Programming}

A mixed-integer linear programming problem (MIP) allows both for continuous and discrete decisions in an optimization problem with linear inputs.  While MIPs are less common in conservation planning, continuous decision variables such as
currency \cite{jafari2013new,jafari2017achieving}, and amount of land or number of parcels allocated to a species \cite{beaudry2016identifying} have been used to decide candidate sites.

 \subsubsection{Mixed Integer Nonlinear Programming}
An optimization problem that simultaneously allows for nonlinear functions, continuous decisions, and discrete decisions is called a mixed integer nonlinear programming problem (MINLP).
\begin{equation}
\begin{aligned}
\label{minlp}
    \underset{x,y}{\textrm{minimize}} \quad &f(x,y) \\[.25em]
    \textrm{subject to} \quad &g(x,y) \leq 0 \\[.25em]
    & x \in X \subseteq \mathbb{Z}^p\\[.25em]
    & y \in Y \subseteq \mathbb{R}^n
\end{aligned}
\end{equation}

where $X$ is the set of discrete variables,  $Y$ is the set of continuous variables, and $f: (X,Y) \mapsto \mathbb{R}$ and $g: (X,Y) \mapsto \mathbb{R}^m$ are sufficiently smooth functions.

One of the first uses of MINLPs in conservation decision-making introduced a nonlinear objective function \cite{ferson2000mathematical} for clustering habitats in a reserve.  The same approach was also used in marine PAs \cite{leslie2003using}. \cite{stralberg2009optimizing} presented an MINLP which had a nonlinear objective and represented salinity in wetland restoration as a continuous variable. While formulations of MINLPs exist, researchers often linearize the nonlinear functions before solving the problem, due to reasons such as availability of modeling environments and solvers, availability of high quality linearizations, and theoretical guarantees for optimality.

\subsubsection{The Proposed MINLP approach} As policymakers push for $30\times 30$, we anticipate a need for more complex spatial optimization models for PAs. With the gap between theory and practitioners in implementing spatial conservation models \cite{sinclair2018use, ferraz2021bridging}, there is a need for interdisciplinary collaboration. Therefore, the goal of this work is to provide a new mathematical programming tool for conservationists which allows for linear and nonlinear inputs, as well as continuous and discrete variables. 

Our proposed MINLP model will be structured with a ``plug-and-play'' feature, where a user would be able to plug in their preferred metric (and the software for computing it) to optimize over. This would allow our model to be used as a framework and provide flexibility to the ongoing advancements. In this paper, we will present a proof-of-concept by connecting the model to a PVA software to represent several nonlinear problem components. Connectivity, competition, crowding, and other similar concerns are handled by the PVA software, rather than expressed as constraints of the optimization model.

\section{Methods}

For this proof-of-concept paper, we will focus only on binary decisions.  The heuristic used in the paper can solve general MINLPs, which will be handled in future research.  In this paper, we will consider two forms of \eqref{minlp}.  The first is a single-objective, constrained problem formulated as follows:
 \begin{align} \label{cons}
    \displaystyle     \underset{x \in \{0,1\}^{|P|}} {\textrm{minimize}} \quad f(x) \quad \textrm{subject to} \quad g(x) \geq \tau 
\end{align}
The second is a multi-objective problem reformulated as an unconstrained problem:
  \begin{align} \label{uncons}
    \displaystyle     \underset{x \in \{0,1\}^{|P|}} {\textrm{minimize}} \quad \sum_{k=1}^{m} \lambda_k f_k(x) 
\end{align}
where $m$ is the number of objectives and $\lambda$ is a vector of penalty parameters.

 In this section, we will discuss these models in further detail as applied to PVA for a metapopulation. The notation can be found in Table \ref{tab:notation}.

\begin{table}
    \centering
\begin{tabular}{lp{0.75\linewidth} }
 \hline
\rule{0pt}{1\normalbaselineskip} 
$P$ & Set of parcels. \\[1em]
$B$ & Landscape divided into $n\times n$ parcels, each with a habitat suitability index, representing the level that a parcel can support this species. Refer to Figure \ref{fig:maps}. \\[.75em]
$X_p$& Binary decision variable denoting whether or not parcel $p \in P$ is preserved. Refer to Figure \ref{fig:maps}.\\[.75em]
$Z$& Resulting habitat from configuration $X$. For any parcel $p \in P$, $Z_p = B_p X_p$.  Refer to Figure \ref{fig:maps}.\\[.75em]
$c_p$	&	Acquisition cost of preserving parcel $p$ \\[.75em]
$r(Z)$ & Risk of total extinction for habitat $Z$.\\[.75em]
$t(Z)$ & Median time to extinction for habitat $Z$.  \\[.75em]
$a(Z)$ & Expected minimum abundance for habitat $Z$.  \\[.75em]
$\rho$ & Threshold by which the solution from $Z$ can differ from the solution given by $B$.\\[.75em]
$\beta$ & Budget\\
 \hline
\end{tabular}
    \caption{Descriptions of the decision variable and parameters used in our models. Refer to the technical appendix for illustrations of functions $r,t,$ and $a$.}
    \label{tab:notation}
\end{table}
For our analysis, we will use the software \emph{RAMAS GIS} \cite{ramas2013_spatial} and \emph{RAMAS Metapopulation} \cite{ramas2013_metapop} as the PVA component. From this software, we opted to use the following metrics in our model: risk of total extinction, time to extinction, and expected minimum abundance. To ensure reproducibility of our results \cite{morrison2016repeatability}, we provide further descriptions of these PVA terms and more detail on implementation can be found in the technical appendix. 

Recall that PVA is based on predictions which we obtain from repeatedly simulating a iteration. In our analysis, we set the duration to be 100 years and generate 1000 iterations for each simulation. The statistics are reported using a 95\% confidence interval based on the Kolmogorov-Smirnov test statistic.  The width of the confidence interval is a function of number of iterations, and to limit the variance in risk we chose the maximum value allowed by the software.

\subsubsection{Median Time to Extinction}
In \emph{RAMAS}, the time to quasi-extinction is the number of years until a metapopulation abundance drops below a given threshold. In our case, we set this threshold as 0 and refer to it as time to extinction.

In each of the 1000 iterations, the time to extinction is recorded. The median is then computed over these iterations. We report this value and refer to it as the median time to extinction.

If a habitat does not reach extinction in at least 500 iterations, a median is not reported. For such cases, the median time to extinction would be recorded as $>100$. That is, in our time frame of 100 years, the population did not go extinct in at least half of the iterations.

\subsubsection{Risk of Extinction}
The terminal extinction risk is the probability that a metapopulation abundance will be below a certain threshold at the end of 100 years. The risk of extinction is the terminal extinction risk when the threshold is chosen to be 0. 

\subsubsection{Expected Minimum Abundance}
For each iteration in the simulation, the abundance of every metapopulation is totaled for each year, and the smallest abundance level over the 100 years is recorded. The expected minimum abundance is the average over 1000 iterations.  This metric is useful to evaluate iterations where the metapopulation may not be close to extinction.

Many conservation models use functional and structural connectivity as constraints. 
For our proposed model, both types of connectivity are considered implicitly through expected minimum abundance, since model parameters that impact abundance (such as carrying capacity, relative survival, relative fecundity) are functions of total habitat suitability (THS).  Larger patches have higher THS, so maximizing or constraining to higher values of expected minimum abundance imposes connectivity. 

We also experimented with average habitat suitability, and found that this encouraged small patches of high habitat suitability that were fragmented. This area could be explored further, as smaller patches are also beneficial and have shown to have higher biodiversity  \cite{fahrig2020several}. 

\begin{figure}
    \centering
\includegraphics{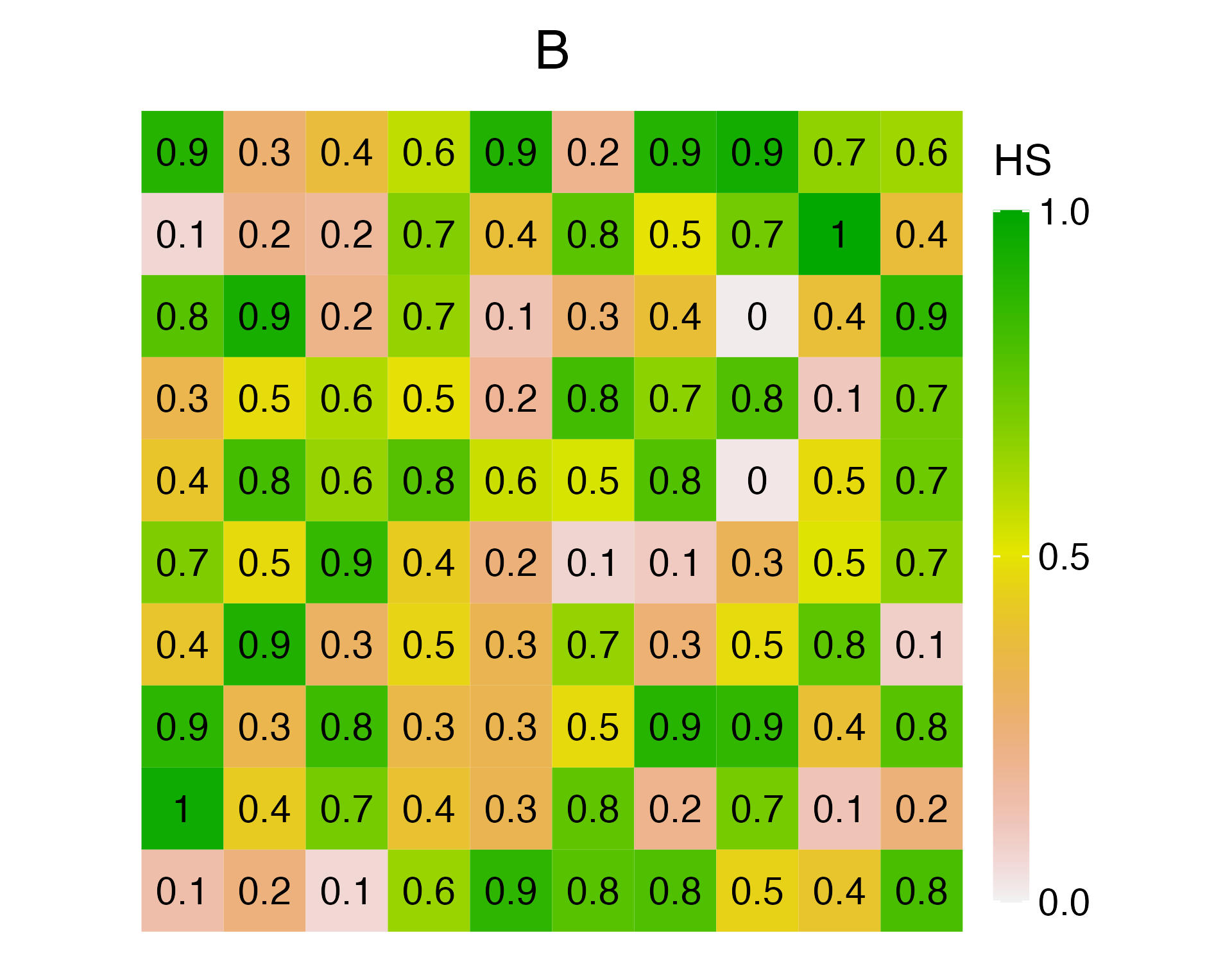}
\includegraphics{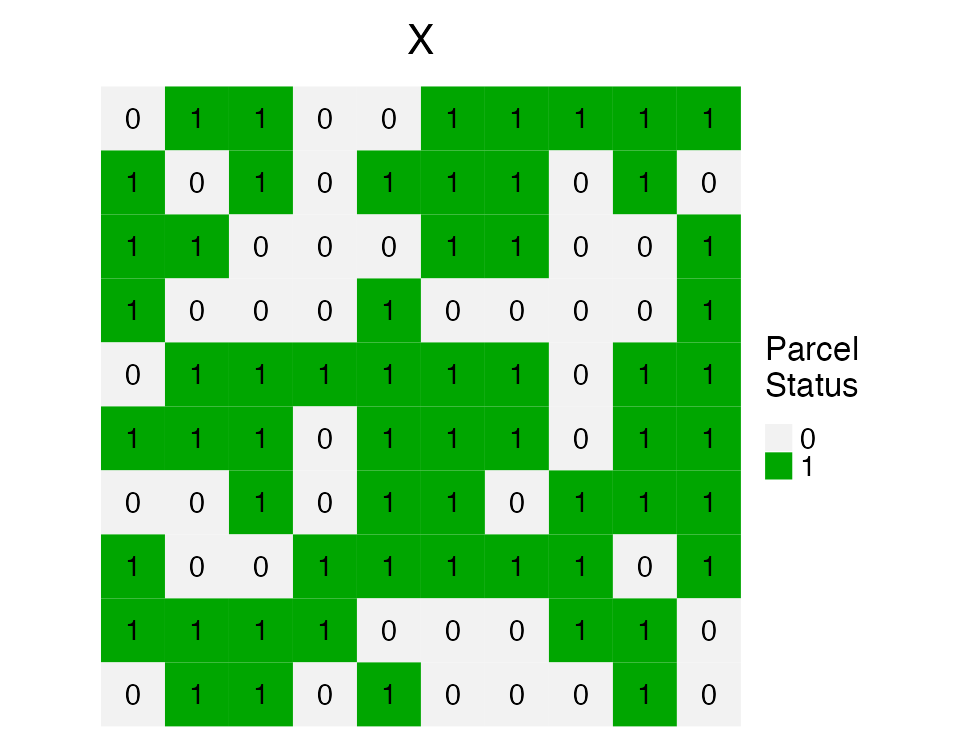}
\includegraphics{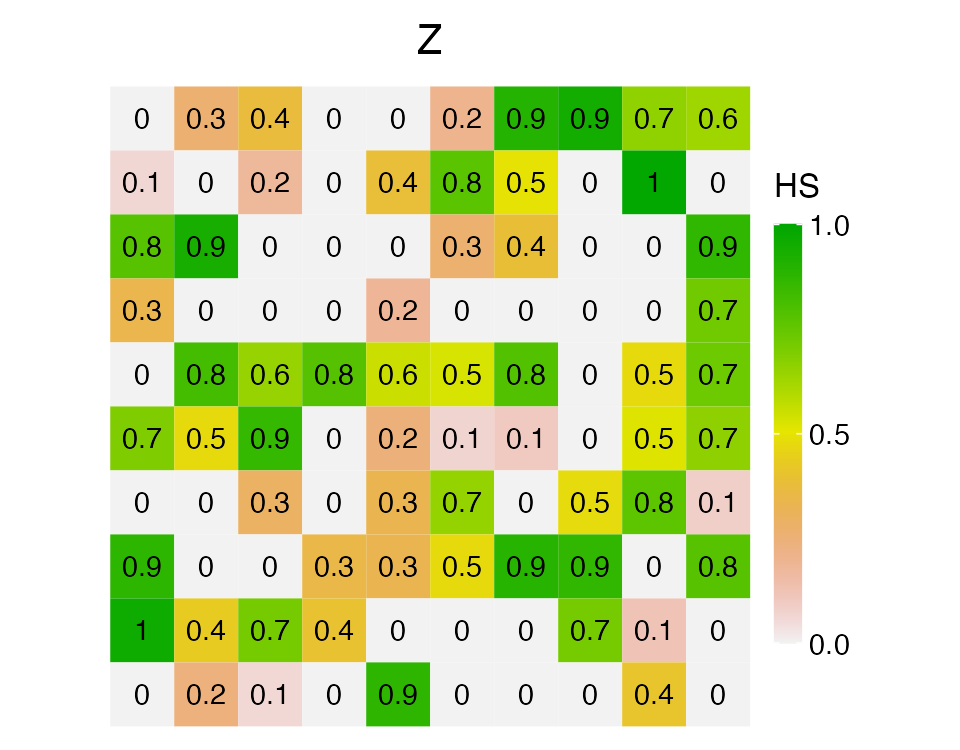}
        \caption{This figure illustrates the process of obtaining the habitat map $Z$. Only the parcels selected by $X$, which is the solution from \eqref{eq:conModel} or \eqref{eq:unconModel}, are retained from the full habitat map $B$. The maps $B$ and $Z$ display the habitat suitability (HS) for each parcel.}
    \label{fig:maps}
\end{figure}

\subsection{Model with Constraints}

The baseline scenario entails preserving every parcel in the landscape $B$. Then, our optimization model with constraints (adapted from \eqref{cons}) aims to obtain the cheapest solution that yields a median time to extinction, risk of extinction, and expected minimum abundance within a pre-defined gap of the corresponding metric for the baseline. The resulting constrained MINLP is as follows:
\begin{equation} \label{eq:conModel}
	\begin{array}{llll}
    \underset{X \in \{0,1\}^{|P|}} {\textrm{minimize}} & \displaystyle \sum_{p \in P} c_p X_p& &  \\[2.1em]
	\textrm{subject to } &  r(Z) - (\rho_r +r(B)) &\leq & 0 \\[1em]
     	 &  t(Z) - \rho_t t(B) &\geq&  0 \\[1em]
    	   &  a(Z) - \rho_B a(B) & \geq  &0 \\[1em]
            & Z_p = B_p X_p    &\forall &p \in P \\[1em]
    	& X_p \in \{0,1\}   &\forall &p \in P
	\end{array}
\end{equation}

The risk of extinction is a probability, so the gap is integrated additively into the model (i.e. risk of extinction for the preserved habitat should be within a small percentage of the baseline risk). Note that in this model, we have a linear objective, discrete decisions, and nonlinear constraints for the PVA metrics obtained from RAMAS, which is invoked as a black-box.

\subsection{Multi-objective Model}

Re-interpreting the cost and the PVA metrics as priorities, we can also formulate a multi-objective model.  Let
\begin{align}
    f(X) = \left[\sum_{p \in P} c_p X_p,r(Z),t(Z), a(Z) \right]
\end{align}
be the vector of objective functions,
and let $\lambda = [\lambda_1,\ldots,\lambda_4]$ be a vector of weights. We use weighted decomposition to convert this $n$ objective problem into a single objective: 
\begin{equation}
    \label{eq:unconModel}
    \begin{array}{ll}
       \underset{X \in \{0,1\}^{|P|}} {\textrm{minimize}} & \displaystyle \sum_{i=1}^4 \lambda_i f_i(X) \\[1.5em]
       {\textrm{subject to }}  & Z_p = B_p X_p  \quad \forall p \in P \\[1em]
       & X_p \in \{0,1\} \quad \forall p \in P
	\end{array}
\end{equation}

\subsection{Solution Heuristic}
The black-box nature of \eqref{eq:conModel} and \eqref{eq:unconModel} yields a MINLP where the underlying nonlinear functions are not well-defined and the overall problems are nonconvex. Therefore, we propose to use a heuristic to solve these problems.  Specifically, we have chosen to use Ant Colony Optimization as our solution method.

\section{Numerical Testing}

\begin{algorithm}[tb]
\caption{Black-box Optimization Framework}
\label{alg:algorithm}
\textbf{Input}: Landscape $B$ \\
\textbf{Parameters}: Model type, $\rho$, $\lambda$\\
\textbf{Output}: $X^*$ 
\begin{algorithmic}[1] 
\IF{Model type is Constrained}
\STATE Obtain baseline PVA metrics $r(B), t(B), a(B)$ for $B$
\STATE Formulate and solve \eqref{eq:conModel} using Ant Colony Optimization, obtain $X^*$
\ELSIF{Model type is Multi-Objective}
\STATE Formulate and solve \eqref{eq:unconModel} using Ant Colony Optimization, obtain $X^*$
\ENDIF
\STATE \textbf{return} $X^*$
\end{algorithmic}
\end{algorithm}

\subsection{Software \& Hardware}
The machine used had the following specifications.
\begin{description}
    \item[Processor] Intel(R) Core(TM) i5-3470S CPU @ 2.90GHz
    \item[Installed RAM] 12 GB 
    \item[OS] Windows 10
    \item[System type] 64-bit OS 
\end{description}
This machine was accessed using a remote desktop on a macOS 13.4.1 with a 3.6 GHz 10-Core Intel Core i9 processor.

The code uses a combination of Python (Version 3.10.9), {\sc R} (Version 4.3.1), \emph{RAMAS GIS: Spatial Data} (Version 6.0), \emph{RAMAS GIS: Habitat Dynamics} (Version 6.0), and \emph{RAMAS Metapopulation} (Version 6.0). The {\sc R} package used for generating the data was \emph{raster} (version 3.6-20). The optimization models were solved using \emph{pygmo} (Version 2.18.0).

\subsection{Data}
The input data $B$ is a randomly generated raster where each pixel has a value that is uniformly distributed from $[0,1]$ to represent the habitat suitability. 

\begin{figure}
    \centering
    \includegraphics{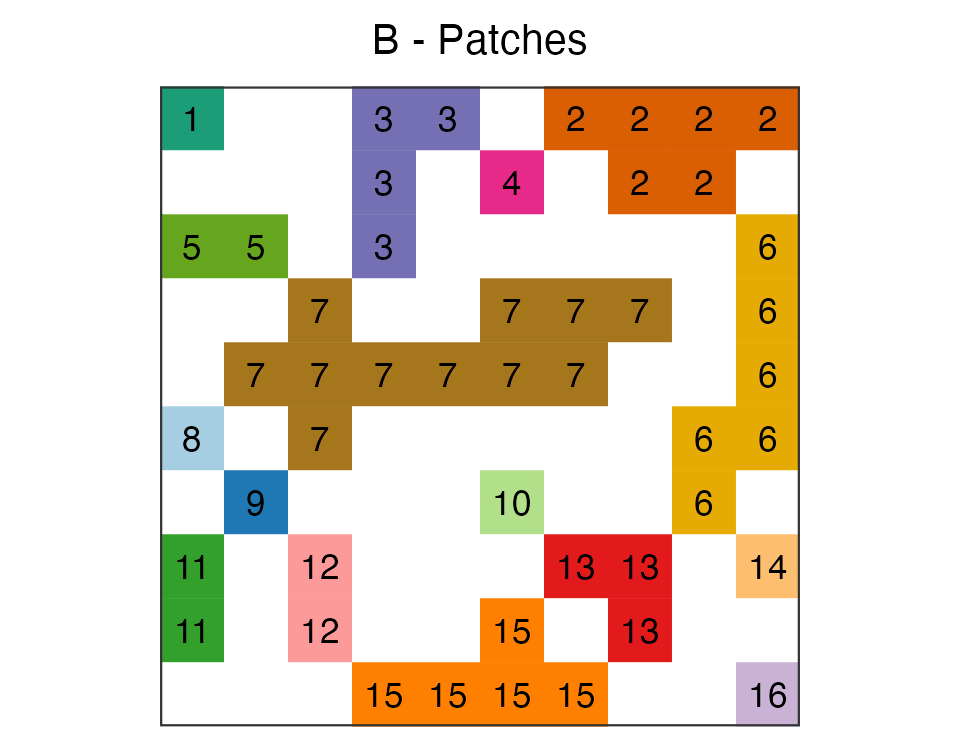}
    \includegraphics{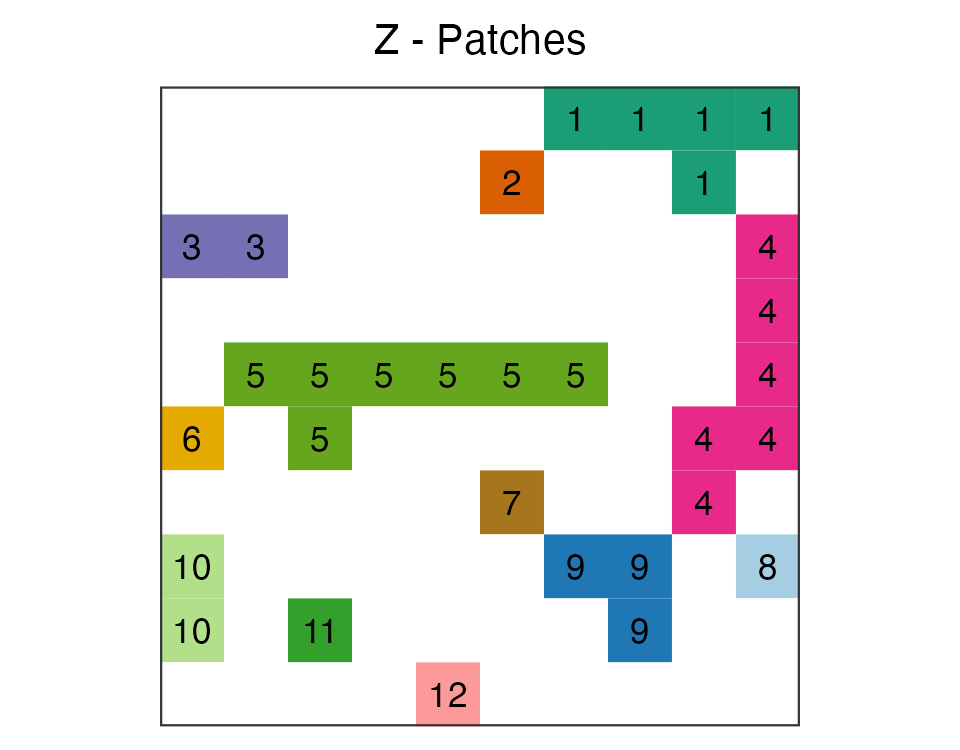}
    \caption{ This figure displays the patch structure of $B$ and $Z$ from Figure \ref{fig:maps}.     
     We see that $Z$ did not preserve patches 1, 3, 9, or 16 from $B$ and it reduced patches 2, 7, 12, and 15. Out of the 48 habitable parcels, 37 are preserved. }
    \label{fig:patches}
\end{figure}

To establish groups of parcels that will contain metapopulations, we set a habitat threshold to $0.5$, where only parcels of higher than the threshold are habitable. A collection of contiguous habitable parcels is denoted as a patch. Refer to Figure \ref{fig:patches}.


\subsection{Testing}

 The following process is described in Algorithm \ref{alg:algorithm} and further detailed in the technical appendix.

Depending on the model type, either \eqref{eq:conModel} or \eqref{eq:unconModel} is solved using Extended Ant Colony Optimization (ACO) as described in \cite{schluter2009extended} and implemented in the Python package, \emph{pygmo}. (The parameters used in ACO are outlined in the technical appendix.) This heuristic is state-of-the-art and it was chosen because the MINLPs are nonconvex nonlinear problems with black-box function evaluations.  In each iteration of the ACO, we create a habitat map $Z$ and the simulation in {\em RAMAS} runs for 100 years and repeats 1000 times to obtain the PVA metrics. 

To create a habitat map $Z$, we begin with the full habitat map $B$ and consider the preserved parcels as indicated by $X$. For parcels that are not preserved, the habitat suitability level decreases such that it is no longer habitable. This is shown in Figure \ref{fig:maps}.

We present two landscape sizes: $10 \times 10$ and $20 \times 20$. For each landscape, we solve \eqref{eq:conModel} and \eqref{eq:unconModel} and report the metrics and those from the baseline scenario in Table \ref{tab:THS_weighted}.

\subsection{Results}

The median time to extinction for each solution was identical to that of the baseline. This indicates that at least half of the simulations did not end in extinction. 

The constrained model for both sizes and the multi-objective model for $n=20$ also have the same risk of extinction as the baseline. Whereas, the risk of extinction has a marginal increase with the multi-objective model for $n=10$. Recall the 95\% confidence interval is based on the Kolmogorov-Smirnov test statistic, thus has a width of $\pm 3\%$ given 1000 replications.
Considering this, a 0.017\% risk of extinction is quite similar to $B$'s risk of 0\%.

For the expected minimum abundance, we begin to see a reflection of patch size on the population level. Habitats with smaller patches will have less abundance. Thus, the habitats with the most abundance are the most expensive. 

The constrained model used $\rho = [0.1, 0.9,0.8]$, thus the constraints are $r(Z) \leq 0.1$, $t(Z) \geq 90.9$, and $a(Z) \geq 61.4$ for $n=10$ and $r(Z) \leq 0.1$, $t(Z) \geq 90.9$, and $a(Z) \geq 271.88$. for $n=20$. ACO does not guarantee feasibility, thus the optimal solution is one which has the minimal constraint violation. We see the risk and time constraint are satisfied, while the abundance is not. 

Comparing the abundance and cost to the multi-objective model, we find that abundance is much higher for the constrained model. This is due to abundance being a constraint, rather than an objective, as satisfying the constraints are prioritized by the solution heuristic.

It follows that the costs of the multi-objective solutions are the cheapest. The weights $\lambda = [0.35,0.15,0.15,0.35]$ valued cost and abundance the same and higher than risk and time. These parameters were selected to reflect the difficulty of obtaining feasible abundance.

The corresponding habitat map for the $10\times 10$ case of $ Z_c^*$ is presented in Figure \ref{fig:maps}, along with the patches in Figure \ref{fig:patches}. Additional figures from Table \ref{tab:THS_weighted} are in the technical appendix.

\begin{table}
    \centering
    \begin{tabular}{llll}
     & $Z_c^*$ & $Z_m^*$ & $B$\\[.1em]
     \hline
    \rule{0pt}{1\normalbaselineskip} 
   $10\times 10$  &  &  & \\[.25em]
    Total Cost & $340$  & $199$ & $577$ \\ 
    Risk & $0$  &  $0.017$ & $0$  \\ 
    Time & $>100$ &$>100$ & $>100$ \\ 
    Abundance  & $33.5$ & $12.1$ &  $61.4$\\ 
    \hline
    \rule{0pt}{1\normalbaselineskip} 
   $20\times20$  &   &  & \\[.25em]
    Total Cost &  1134 & 874 & 2137\\ 
    Risk  & 0  &  0 & 0 \\ 
    Time  & $>100$  &  $>100$& $>100$\\ 
    Abundance   &  135.2 & 69.4 & 339.6 \\ 
    \end{tabular}
    \caption{Let $Z_c^*$ and $Z_m^*$ be the landscape given from $X^*$ in \eqref{eq:conModel} and \eqref{eq:unconModel}, respectively. $B$ is the case if every parcel were preserved. \emph{Risk} is the risk to total extinction. \emph{Time} is the median time to extinction. \emph{Abundance} is the expected minimum abundance.}
    \label{tab:THS_weighted}
\end{table}

\section{Discussion}

\subsubsection{Limitations}

Due to our selected PVA software, we were unable to solve using parallel computing, which drastically impacted the runtime. For this reason, we needed to limit the size of the landscape and the number of ACO generations, as these have the biggest impact on runtime. 

Let $g$ denote the number of generations, the number of iterations are then 
\[\textrm{Iterations} = (n^2+1)*(g+1)\]

Without parallel computation, the problem scales by a factor of $n^2$. This is significant considering that each iteration requires around 31 and 36 seconds for $n=10$ and $n=20$, respectively. 

Ideally, the PVA component (or other selected software) would need to allow parallelization when working with larger landscape sizes. 
While long-term conservation decisions do not necessarily need to be made in a short window, parallelization can speed up the numerical testing and model refinement process.

Additionally, the MINLP with black-box functions limits us to using solution heuristics, which lack theoretical optimality guarantees. 

Still, heuristics can offer computational advantages and are not unusual in conservation settings.
MARXAN \cite{ball2009marxan}, a population conservation planning software, uses simulated annealing. 
Even with explicitly defined convex formulations in optimization models, heuristics can be preferred for their interpretability and ability to scale up through parallelization.

Lastly, this framework is intended as a starting point to further multidisciplinary research. AI should aid conservationists, not replace them.  The tools that are helpful for decision-makers must have the ability to tackle complex scenarios and integrate their software and models seamlessly.

\subsubsection{Future Directions }
There are a few avenues to explore in future work. We have already indicated that parallelization of the optimization method will be a priority. We can also incorporate methods from derivate-free optimization literature, such as using surrogate functions for the PVA metrics.

The testing in this paper used randomly generated data. While this does not impact the mechanisms of the framework, the overall goal is to use real data that reflect real-world problems. 
Species distribution models produce range maps depicting the geographic locations where a species may occur and are often used in spatial optimization models. Applying such methods would provide more validity to the results.

Finally, another area of exploration would be incorporating metrics that allow for multiple species. This implementation would increase problem size, as every species would add another dimension to the model.  It would also increase opportunities and motivation for parallelization of the optimization method.

\section*{Ethical Statement}

Historically, protected areas have been used as a tool for colonization under the guise of conservation, largely impacting the Indigenous Peoples \cite{west2006parks,stevens2014indigenous}.  
 Optimization models in conservation planning have potential to yield inequitable decisions, particularly in the context of transnational conservation initiatives such as $30\times 30$ \cite{chapman2021promoting}, 
thus it is imperative that this framework--and others like it--are used with discretion. 

\section*{Acknowledgements}




We would like to thank the anonymous reviewers for their valuable feedback.
This work was supported by the LeBow College of Business Innovation Micro-Grant for DEI \& Environmental Sustainability at Drexel University.

\bibliography{ref.bib}

\begin{appendices}

  \input{appendix}

\end{appendices}

\end{document}

%% file: appendix.tex
\section{Technical Appendix}
\subsection{Data \& Code Availability}
\begin{figure}[p]
    \centering
\includegraphics[width = 0.35\textwidth]{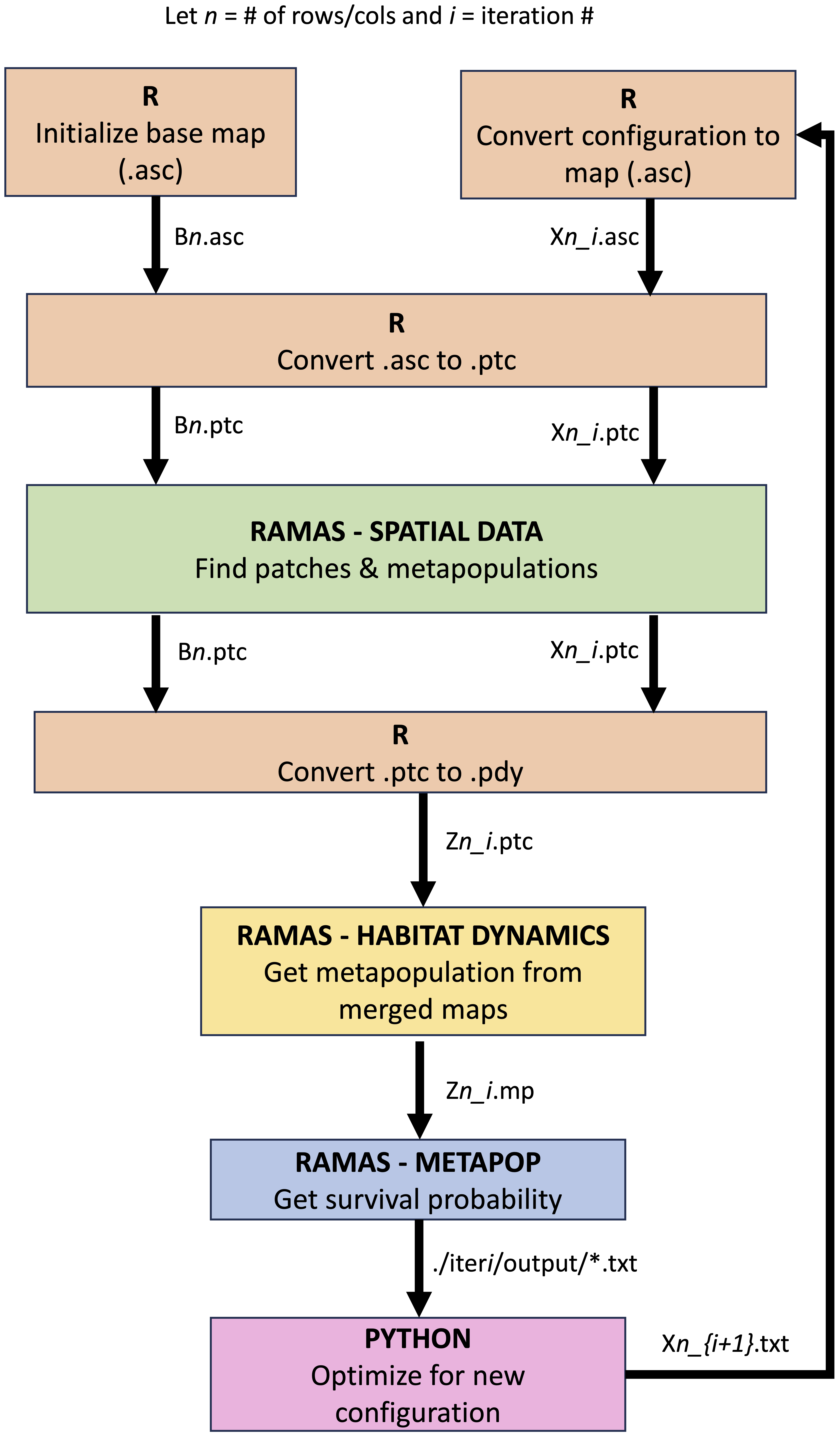}
    \caption{Our software using PVA broken down by each step in the code.}
    \label{fig:flow}
\end{figure}

The data and code can be found at the following repository https://github.com/cassiebuhler/conservation-dfo. For an illustration of the file formats, see Figure \ref{fig:flow}.



\subsection{Framework Parameters}

\subsubsection{Ant Colony Optimization}
The parameters for our \emph{pygmo} optimization solver, ACO \cite{schluter2009extended}, are in Table \ref{tab:ACO}.

\begin{table}[ht]
    \centering
    \begin{tabular}{ll}
        \hline
    \rule{0pt}{1\normalbaselineskip}\noindent
      \textbf{Ant Colony Optimization} & \\
        \hline
    \rule{0pt}{1\normalbaselineskip}\noindent
    Generations & 4 \\ 
    Kernel & $n$\\
    Impstop & 1000\\
    Oracle & 1e6
    \end{tabular}
    \caption[ACO parameters in the MINLP framework]{\emph{Generations} is the number of generations to run the algorithm, \emph{Kernel} is the number of solutions stored in archive, \emph{Impstop} stops the algorithm when this number of runs occurs without improvements, and \emph{Oracle} is the penalty parameter. }
    \label{tab:ACO}
\end{table}

\subsubsection{PVA Software}

The PVA component of our framework was implemented with \emph{RAMAS GIS} \cite{ramas2013_spatial} and \emph{RAMAS Metapopulation} \cite{ramas2013_metapop}. In \emph{RAMAS GIS}, we used \emph{RAMAS GIS: Spatial Data} for generating our habitat and metapopulation, and \emph{RAMAS GIS: Habitat Dynamics} to simulate the effects of parcels being protected. The parameters for \emph{RAMAS GIS: Spatial Data} are in Table \ref{tab:spatial1} and \emph{RAMAS GIS: Habitat Dynamics} in Table \ref{tab:hab}. They are separated into their respective categories by where the input appears in the RAMAS interface. 

Then, PVA was conducted in \emph{RAMAS Metapopulation}, where each habitat configuration was evaluated for 100 years, and replicated 1000 times.


\begin{table}[h]
\centering
    \begin{tabular}{ll}
        \hline
    \rule{0pt}{1\normalbaselineskip}\noindent
 \textbf{RAMAS GIS: Habitat Dynamics} & \\  
         \hline
    \rule{0pt}{1\normalbaselineskip}\noindent 
    \textbf{Habitat Changes} & \\
       B - Time step & 1\\
        B - K change to next & same until next \\
        B - F change to next & same until next \\
        B - S change to next & same until next \\
       Z - Time step & 10\\
    \end{tabular}
    \caption[RAMAS GIS: Habitat Dynamics parameters in the MINLP framework]{Parameters for RAMAS GIS: Habitat Dynamics. We assume the habitat is B for the first 10 years, then use Z for the remaining 90 years. \emph{K} is the carrying capacity, \emph{S} is the relative survival rate, and \emph{F} is relative fecundity rate. These values are dependent on \emph{ths}, thus abundance will change from $B$ to $Z$.   }
    \label{tab:hab}
\end{table}
\subsubsection{Density Dependency}
Scramble is defined in the RAMAS User Manual (6.0) as ``Logistic or Ricker type of density dependence, characteristic of the effect of scramble competition.''

The deterministic growth rate based on population size at time step $t$ is the following
\begin{align}
    R(t) = R_{\textrm{max}} * \exp{\left(\frac{\ln(R_{\textrm{max}})*N_p}{K(t)}\right)}
\end{align}
where $N_p$ is the current total number of individuals in a population $p$, and $K(t)$ is the carrying capacity at time $t$. Because our data is randomly generated, this type was also randomly selected.

\subsubsection{Dispersal and Correlation}
The dispersal rate (migration) between population $i$ and population $j$ is calculated as
\begin{align}
    a * \exp{\left(\dfrac{-D_{ij}^{c}}{b} \right)}.
    \label{eq:dispcorr}
\end{align}
This function is also used as the correlation of population fluctuation between population $i$ and population $j$.


\begin{table}[ht]
    \centering
    \begin{tabular}{ll}
            \hline
    \rule{0pt}{1\normalbaselineskip}\noindent
     \textbf{RAMAS GIS: Spatial Data} & \\   
        \hline
    \rule{0pt}{1\normalbaselineskip}\noindent
    \textbf{Habitat Relationships} & \\
    HS function & [<ASC file>]  \\
    HS threshold & 0.500\\
    Neighborhood distance (cells) & 1.00\\
        \hline
    \rule{0pt}{1\normalbaselineskip}\noindent
    \textbf{Link to Metapopulation} & \\
    Carrying capacity (K) & $ths*4$ \\
       $R_{\textrm{max}}$ & $1.5$\\
       Initial abundance  & $ths*2$ \\
       Relative survival  & $\max(1,ths*1.2)$\\
       Relative fecundity   & $\max(1,ths*1.2)$\\
       Distances & Edge to edge \\
            \hline
    \rule{0pt}{1\normalbaselineskip}\noindent
    \textbf{Default Population} & \\
        Local threshold & 0 \\
        Density dependency type & Scramble \\
        Standard deviation of K & 0 \\
        \hline
    \rule{0pt}{1\normalbaselineskip}\noindent
    \textbf{Dispersal} & \\
        $a$ & 0.5\\
       $b$ & 0.8\\
        $c$ & 1\\
       $D$ & 1\\
        \hline
    \rule{0pt}{1\normalbaselineskip}\noindent
      \textbf{Correlation} & \\
       $a$ & 0.8\\
       $b$ & 2\\
       $c$ & 1 \\ \hline
    \end{tabular}
    \caption[RAMAS GIS: Spatial Data parameters in the MINLP framework]{Parameters for RAMAS GIS: Spatial Data. \emph{ths} is total habitat suitability. The formula for dispersal and correlation is \eqref{eq:dispcorr}}
    \label{tab:spatial1}
\end{table}

\subsection{Additional Figures}
\subsubsection{Visualizing PVA Metrics}

\begin{figure}[p]
    \centering
    \includegraphics[width = 0.35\textwidth]{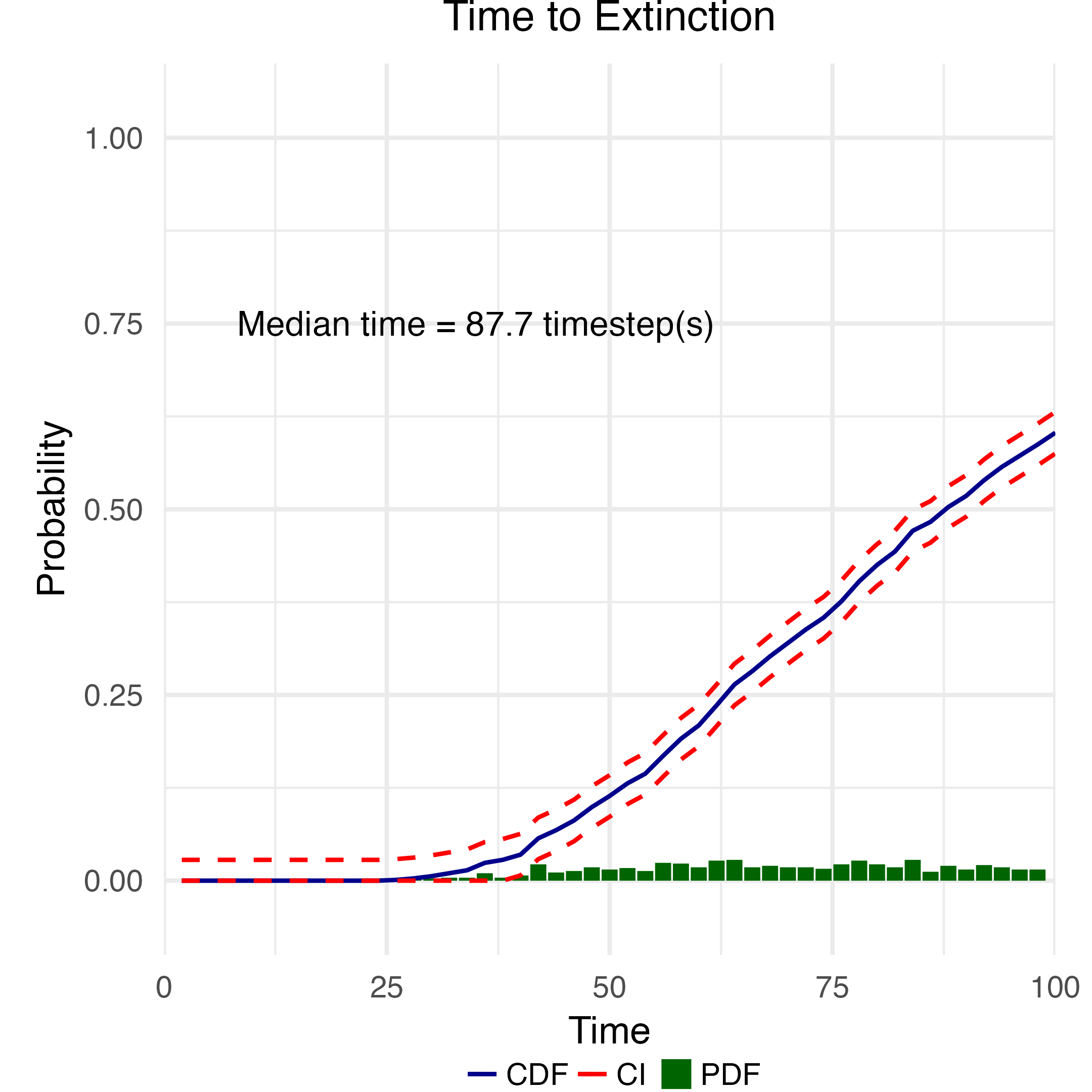} 

    \vspace{0.1in}
    
    \includegraphics[width = 0.35\textwidth]{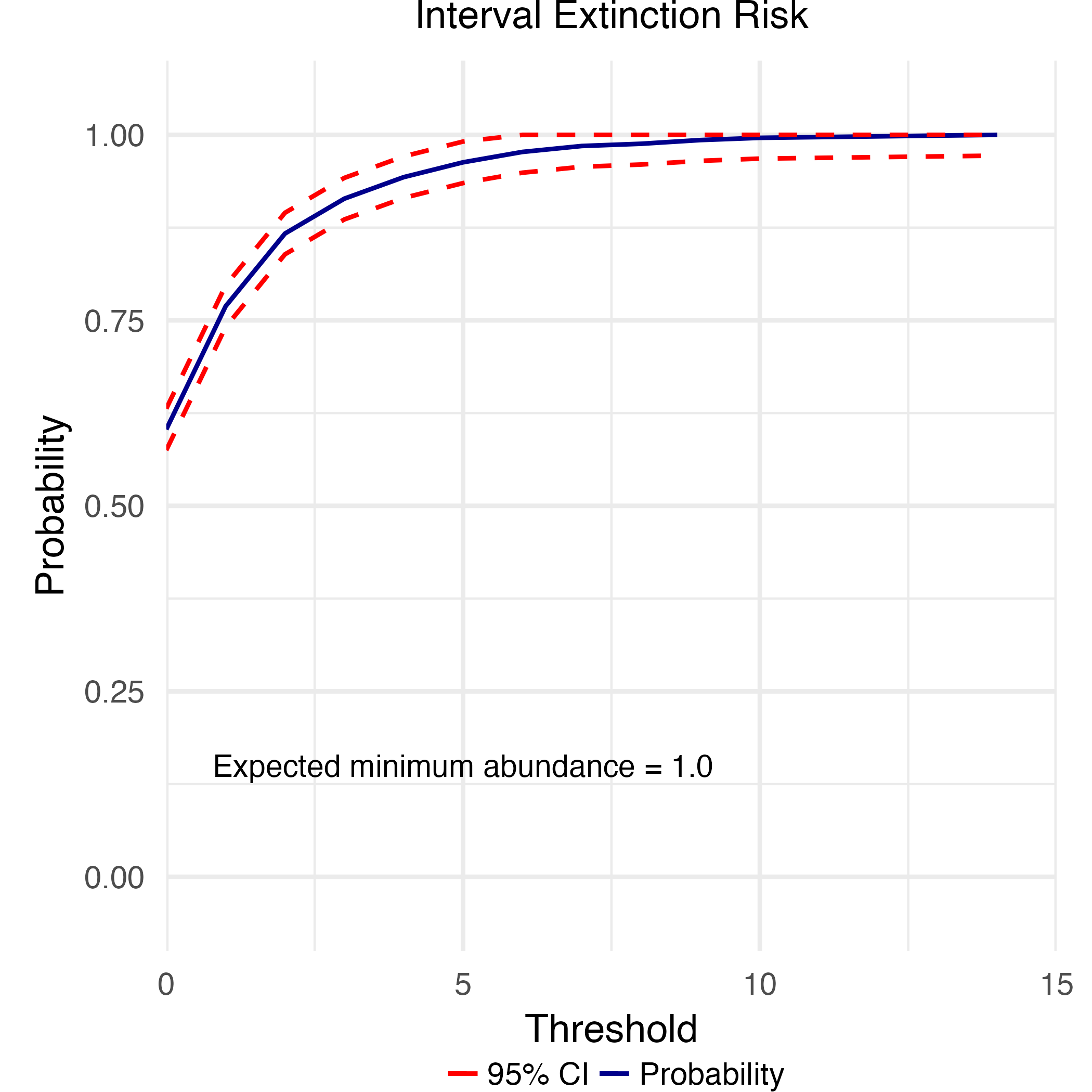}

    \vspace{0.1in}
    
    \includegraphics[width = 0.35\textwidth]{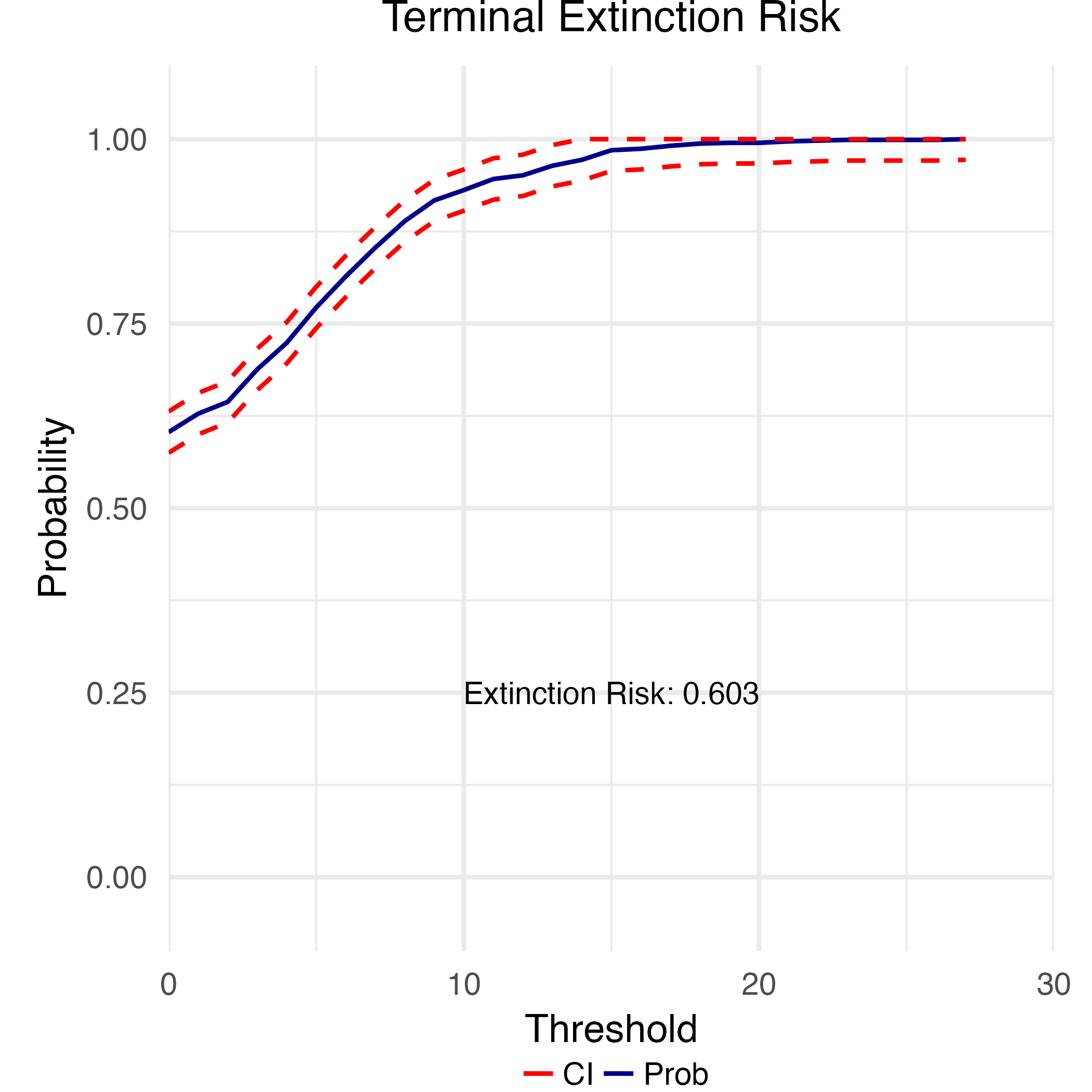}
    
    \caption{Median time to quasi extinction is where the probability is 0.5.  Expected minimum abundance is where the probability is 0.5.  Risk of extinction is the probability when the threshold is 0.}
    \label{fig:n21}
\end{figure}

In Figure \ref{fig:n21}, we visualize PVA metrics with an example of a habitat configuration that was predicted to go to extinction.

\subsubsection{Framework Solutions}
\begin{figure}[p]
    \centering
    \includegraphics{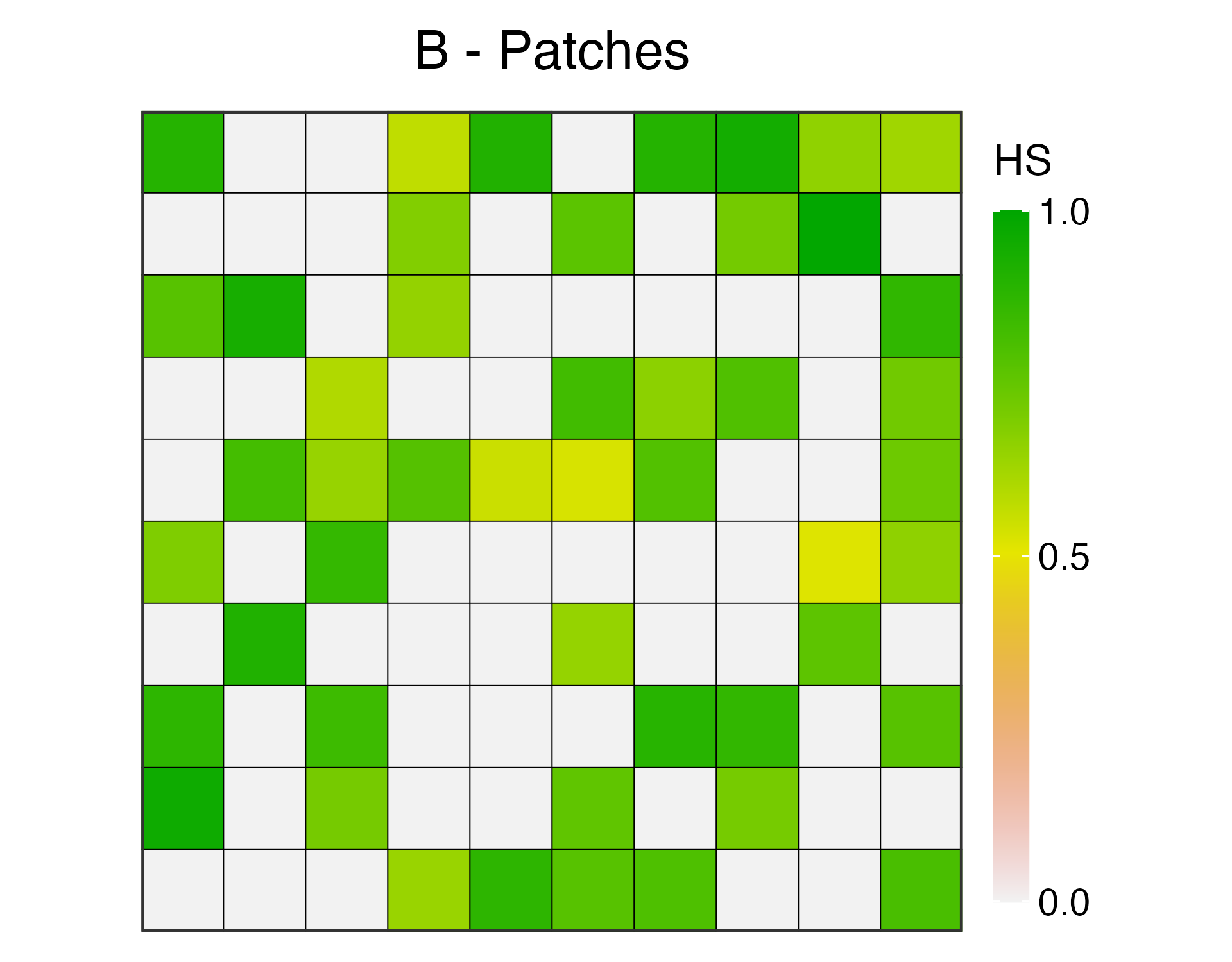}
        \includegraphics{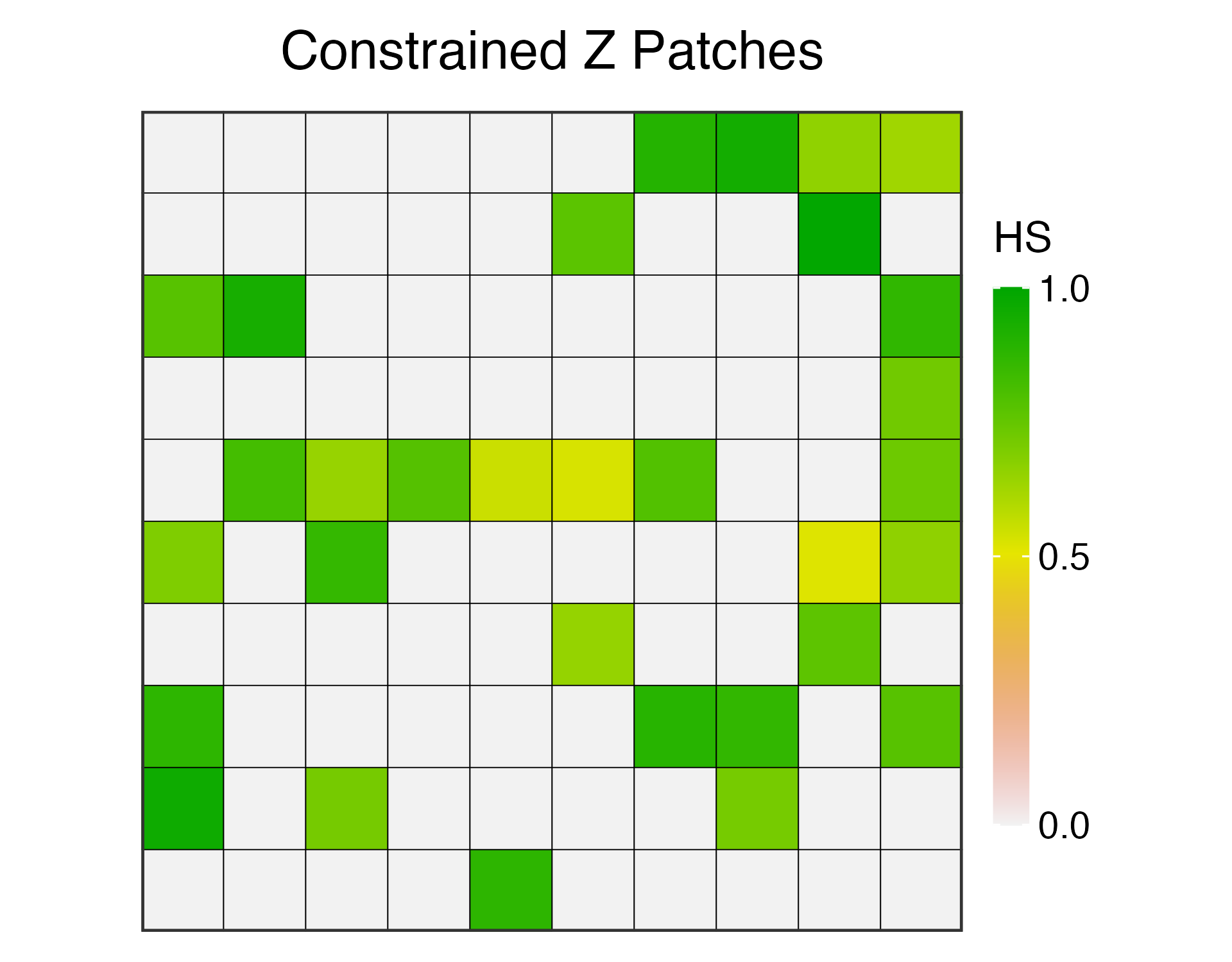}
    \includegraphics{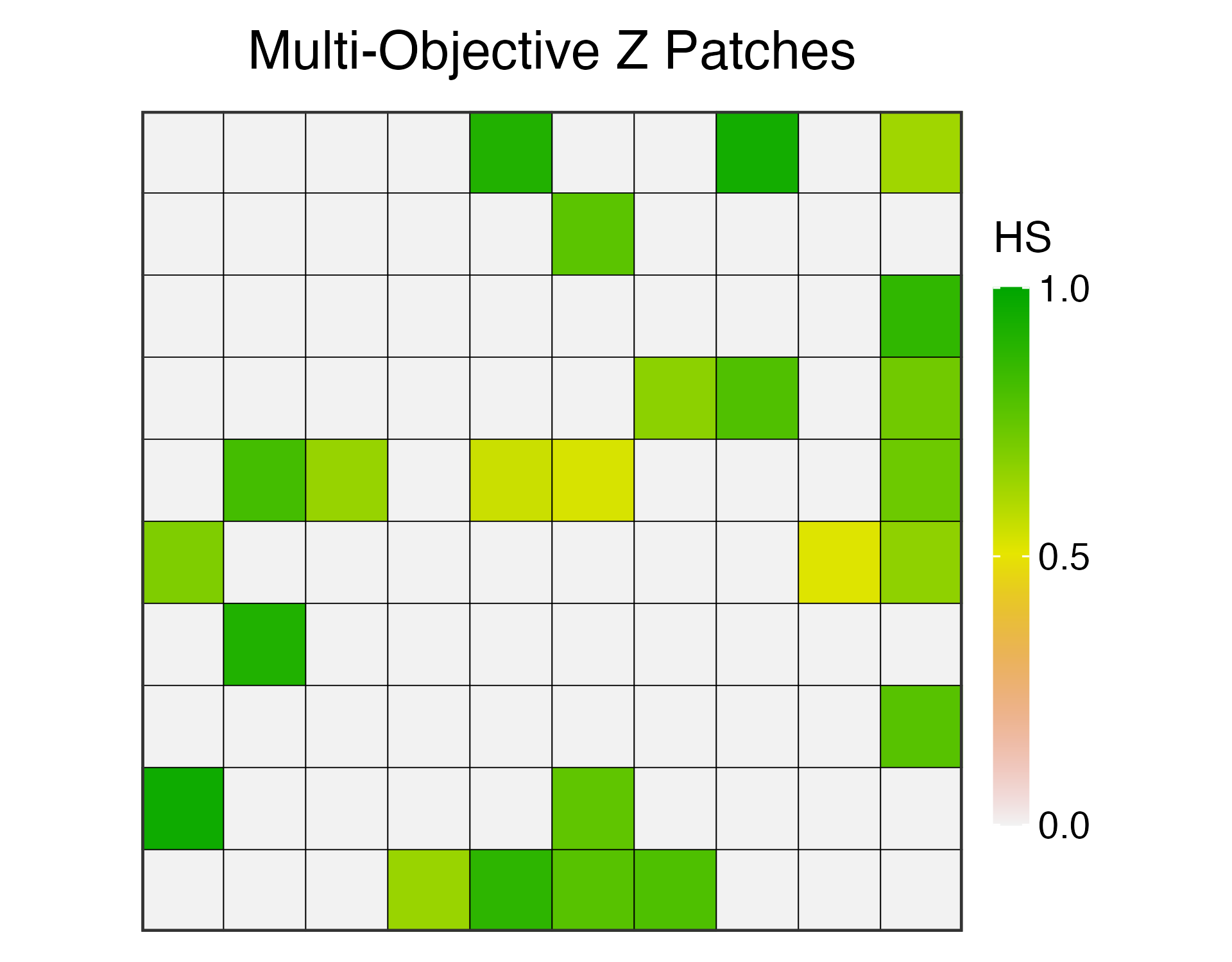}
    \caption{$10\times 10$ solutions from paper. The baseline habitat has 16 populations. The constrained model has 12 populations and the multi-objective model has 13 populations.}
    \label{fig:n10}
\end{figure}

\begin{figure}[p]
    \centering
    \includegraphics{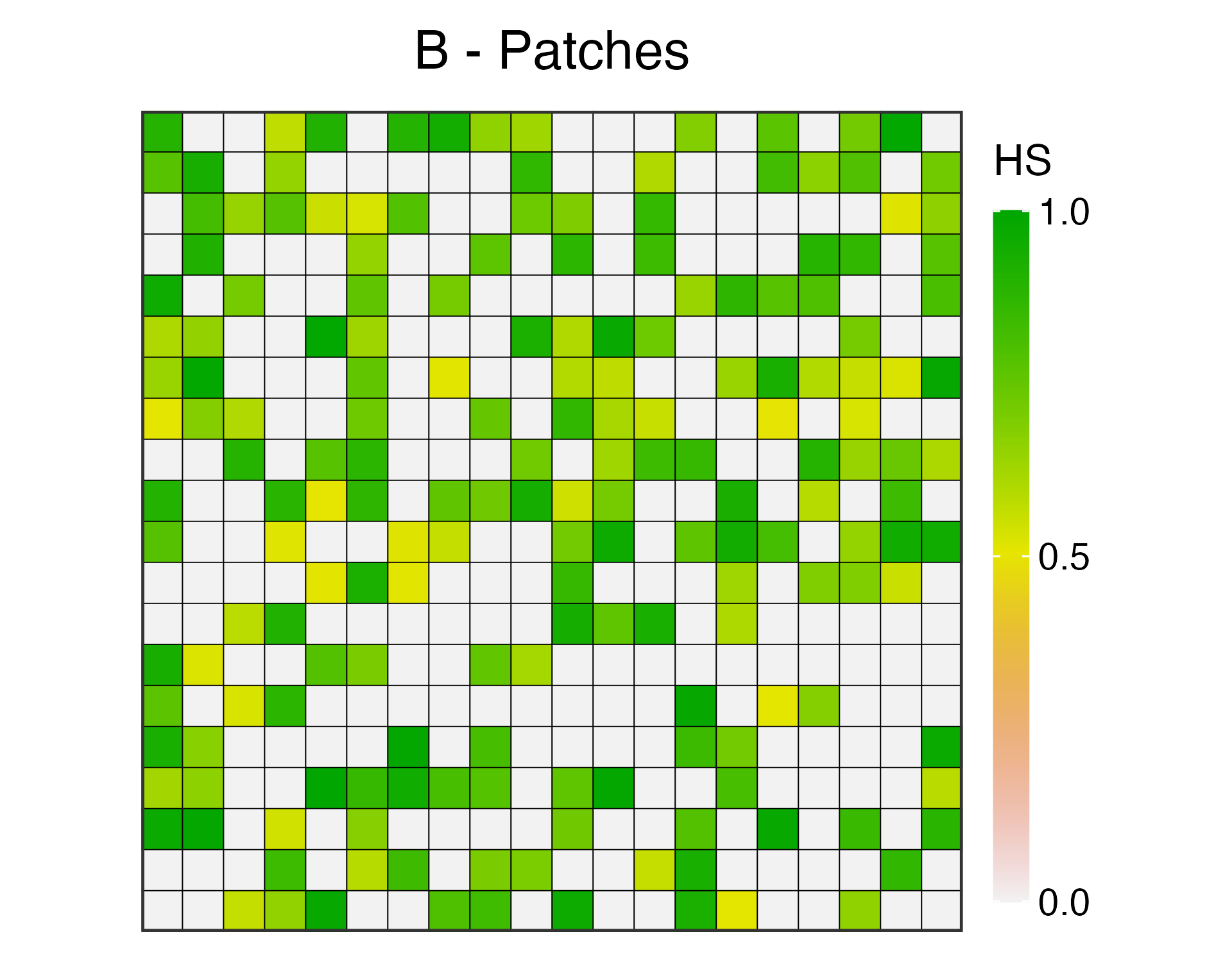}
    \includegraphics{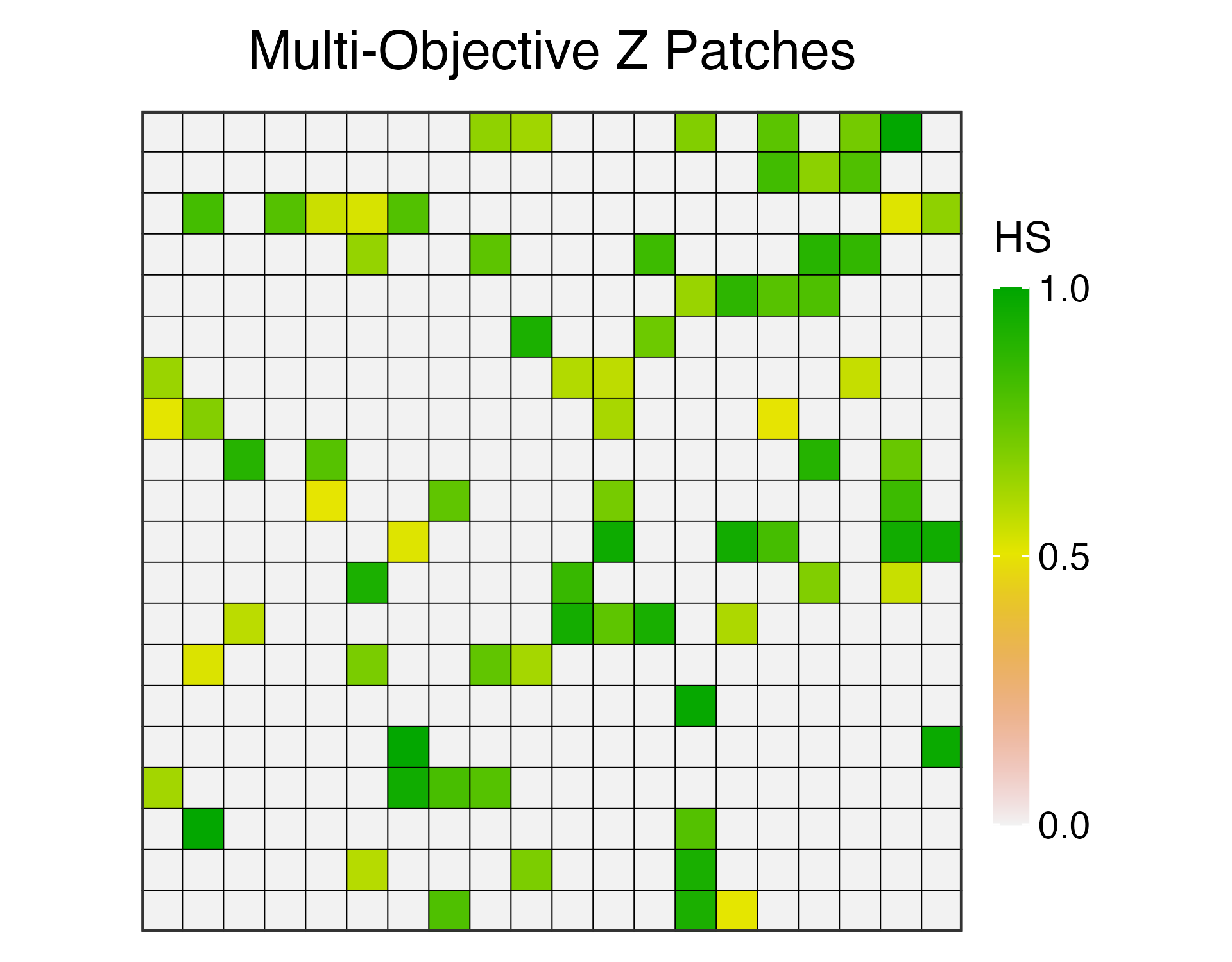}
    \includegraphics{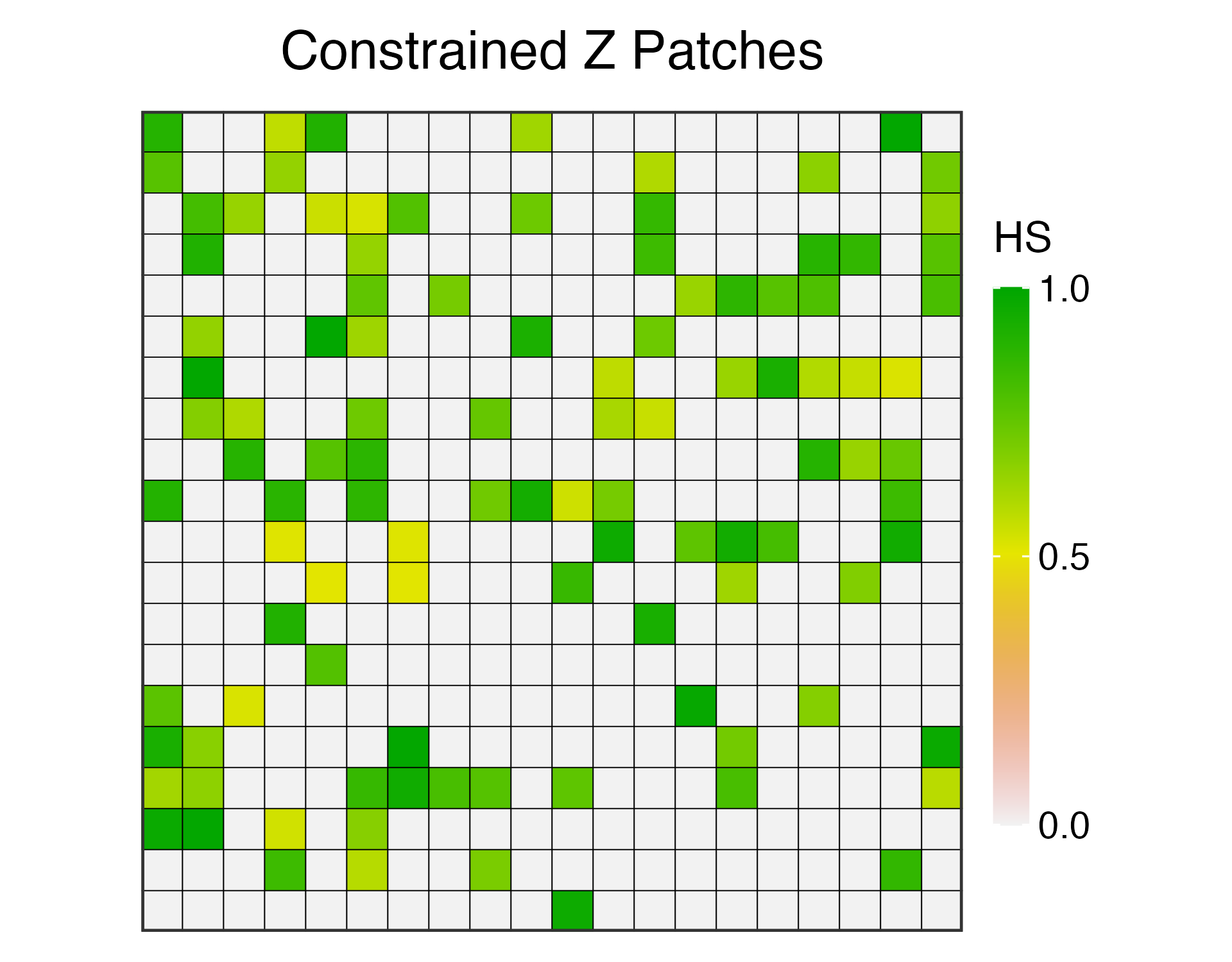}
    \caption{$20\times 20$ solutions from paper. The baseline habitat has 35 populations. The constrained model has 43 populations and the multi-objective model has 48 populations.}
    \label{fig:n20}
\end{figure}

In this paper, we presented the figures for the solution to the $10\times 10 $ constrained model. Here we include the $Z$ configuration, except only where it is habitable. In Figure \ref{fig:n10} we display the $10\times10$ results and in Figure \ref{fig:n20} has the $20\times 20$ results.




%% file: main.bbl
\begin{thebibliography}{59}
\providecommand{\natexlab}[1]{#1}

\bibitem[{Akasaka et~al.(2017)Akasaka, Kadoya, Ishihama, Fujita, and Fuller}]{akasaka2017smart}
Akasaka, M.; Kadoya, T.; Ishihama, F.; Fujita, T.; and Fuller, R.~A. 2017.
\newblock Smart protected area placement decelerates biodiversity loss: a representation-extinction feedback leads rare species to extinction.
\newblock \emph{Conservation Letters}, 10(5): 539--546.

\bibitem[{Ak{\c{c}}akaya and Sj{\"o}gren-Gulve(2000)}]{akccakaya2000population}
Ak{\c{c}}akaya, H.~R.; and Sj{\"o}gren-Gulve, P. 2000.
\newblock Population viability analyses in conservation planning: an overview.
\newblock \emph{Ecological bulletins}, 9--21.

\bibitem[{Ak\k{c}akaya and Root(2013{\natexlab{a}})}]{ramas2013_spatial}
Ak\k{c}akaya, H.~R.; and Root, W. 2013{\natexlab{a}}.
\newblock {RAMAS GIS}: Linking spatial data with population viability analysis.
\newblock \emph{Applied Biomathematics, Setauket, NY}.

\bibitem[{Ak\k{c}akaya and Root(2013{\natexlab{b}})}]{ramas2013_metapop}
Ak\k{c}akaya, H.~R.; and Root, W. 2013{\natexlab{b}}.
\newblock {RAMAS} Metapop: Viability analysis for stage-structured metapopulations (version 6).
\newblock \emph{Applied Biomathematics, Setauket, NY}.

\bibitem[{Alagador and Cerdeira(2022)}]{alagador2022operations}
Alagador, D.; and Cerdeira, J.~O. 2022.
\newblock Operations research applicability in spatial conservation planning.
\newblock \emph{Journal of Environmental Management}, 315: 115172.

\bibitem[{Andersen, Jang, and Borz{\'e}e(2023)}]{andersen2023influence}
Andersen, D.; Jang, Y.; and Borz{\'e}e, A. 2023.
\newblock Influence of landscape and connectivity on anuran conservation: population viability analyses to designate protected areas.
\newblock \emph{Animal Conservation}, 26(3): 381--397.

\bibitem[{Arthur et~al.(1997)Arthur, Hachey, Sahr, Huso, and Kiester}]{arthur1997finding}
Arthur, J.~L.; Hachey, M.; Sahr, K.; Huso, M.; and Kiester, A. 1997.
\newblock Finding all optimal solutions to the reserve site selection problem: formulation and computational analysis.
\newblock \emph{Environmental and Ecological Statistics}, 4(2): 153--165.

\bibitem[{Ball, Possingham, and Watts(2009)}]{ball2009marxan}
Ball, I.~R.; Possingham, H.~P.; and Watts, M. 2009.
\newblock Marxan and relatives: software for spatial conservation prioritisation.
\newblock \emph{Spatial conservation prioritisation: Quantitative methods and computational tools}, 185--195.

\bibitem[{Beaudry et~al.(2016)Beaudry, Ferris, Pidgeon, and Radeloff}]{beaudry2016identifying}
Beaudry, F.; Ferris, M.~C.; Pidgeon, A.~M.; and Radeloff, V.~C. 2016.
\newblock Identifying areas of optimal multispecies conservation value by accounting for incompatibilities between species.
\newblock \emph{Ecological modelling}, 332: 74--82.

\bibitem[{Billionnet(2013)}]{billionnet2013mathematical}
Billionnet, A. 2013.
\newblock Mathematical optimization ideas for biodiversity conservation.
\newblock \emph{European Journal of Operational Research}, 231(3): 514--534.

\bibitem[{Camm et~al.(1996)Camm, Polasky, Solow, and Csuti}]{camm1996note}
Camm, J.~D.; Polasky, S.; Solow, A.; and Csuti, B. 1996.
\newblock A note on optimal algorithms for reserve site selection.
\newblock \emph{Biological Conservation}, 78(3): 353--355.

\bibitem[{Chapman et~al.(2021)Chapman, Oestreich, Frawley, Boettiger, Diver, Santos, Scoville, Armstrong, Blondin, Chand et~al.}]{chapman2021promoting}
Chapman, M.~S.; Oestreich, W.~K.; Frawley, T.~H.; Boettiger, C.; Diver, S.; Santos, B.~S.; Scoville, C.; Armstrong, K.; Blondin, H.; Chand, K.; et~al. 2021.
\newblock Promoting equity in the use of algorithms for high-seas conservation.
\newblock \emph{One Earth}, 4(6): 790--794.

\bibitem[{Chaudhary and Oli(2020)}]{chaudhary2020critical}
Chaudhary, V.; and Oli, M.~K. 2020.
\newblock A critical appraisal of population viability analysis.
\newblock \emph{Conservation Biology}, 34(1): 26--40.

\bibitem[{Church, Stoms, and Davis(1996)}]{church1996reserve}
Church, R.~L.; Stoms, D.~M.; and Davis, F.~W. 1996.
\newblock Reserve selection as a maximal covering location problem.
\newblock \emph{Biological conservation}, 76(2): 105--112.

\bibitem[{Coad et~al.(2019)Coad, Watson, Geldmann, Burgess, Leverington, Hockings, Knights, and Di~Marco}]{coad2019widespread}
Coad, L.; Watson, J.~E.; Geldmann, J.; Burgess, N.~D.; Leverington, F.; Hockings, M.; Knights, K.; and Di~Marco, M. 2019.
\newblock Widespread shortfalls in protected area resourcing undermine efforts to conserve biodiversity.
\newblock \emph{Frontiers in Ecology and the Environment}, 17(5): 259--264.

\bibitem[{Cocks and Baird(1989)}]{cocks1989using}
Cocks, K.; and Baird, I.~A. 1989.
\newblock Using mathematical programming to address the multiple reserve selection problem: an example from the Eyre Peninsula, South Australia.
\newblock \emph{Biological Conservation}, 49(2): 113--130.

\bibitem[{{Convention on Biological Diversity}(2020)}]{GDF2020}
{Convention on Biological Diversity}. 2020.
\newblock Update of the zero draft of the post-2020 global biodiversity framework.
\newblock \url{https://www.cbd.int/doc/c/efb0/1f84/a892b98d2982a829962b6371/wg2020-02-03-en.pdf}.

\bibitem[{Costanza et~al.(2020)Costanza, Watling, Sutherland, Belyea, Dilkina, Cayton, Bucklin, Roma{\~n}ach, and Haddad}]{costanza2020preserving}
Costanza, J.~K.; Watling, J.; Sutherland, R.; Belyea, C.; Dilkina, B.; Cayton, H.; Bucklin, D.; Roma{\~n}ach, S.~S.; and Haddad, N.~M. 2020.
\newblock Preserving connectivity under climate and land-use change: No one-size-fits-all approach for focal species in similar habitats.
\newblock \emph{Biological Conservation}, 248: 108678.

\bibitem[{Cruz, Santulli-Sanzo, and Ceballos(2021)}]{cruz2021global}
Cruz, C.; Santulli-Sanzo, G.; and Ceballos, G. 2021.
\newblock Global patterns of raptor distribution and protected areas optimal selection to reduce the extinction crises.
\newblock \emph{Proceedings of the National Academy of Sciences}, 118(37): e2018203118.

\bibitem[{Csuti et~al.(1997)Csuti, Polasky, Williams, Pressey, Camm, Kershaw, Kiester, Downs, Hamilton, Huso et~al.}]{csuti1997comparison}
Csuti, B.; Polasky, S.; Williams, P.~H.; Pressey, R.~L.; Camm, J.~D.; Kershaw, M.; Kiester, A.~R.; Downs, B.; Hamilton, R.; Huso, M.; et~al. 1997.
\newblock A comparison of reserve selection algorithms using data on terrestrial vertebrates in Oregon.
\newblock \emph{Biological Conservation}, 80(1): 83--97.

\bibitem[{Dilkina et~al.(2017)Dilkina, Houtman, Gomes, Montgomery, McKelvey, Kendall, Graves, Bernstein, and Schwartz}]{dilkina2017trade}
Dilkina, B.; Houtman, R.; Gomes, C.~P.; Montgomery, C.~A.; McKelvey, K.~S.; Kendall, K.; Graves, T.~A.; Bernstein, R.; and Schwartz, M.~K. 2017.
\newblock Trade-offs and efficiencies in optimal budget-constrained multispecies corridor networks.
\newblock \emph{Conservation Biology}, 31(1): 192--202.

\bibitem[{Doak et~al.(2015)Doak, Himes~Boor, Bakker, Morris, Louthan, Morrison, Stanley, and Crowder}]{doak2015recommendations}
Doak, D.~F.; Himes~Boor, G.~K.; Bakker, V.~J.; Morris, W.~F.; Louthan, A.; Morrison, S.~A.; Stanley, A.; and Crowder, L.~B. 2015.
\newblock Recommendations for improving recovery criteria under the US Endangered Species Act.
\newblock \emph{BioScience}, 65(2): 189--199.

\bibitem[{Dudley(2008)}]{dudley2008guidelines}
Dudley, N. 2008.
\newblock \emph{Guidelines for applying protected area management categories}.
\newblock {IUCN}.

\bibitem[{{Executive Order No. 14,008}(2021)}]{EO_30}
{Executive Order No. 14,008}. 2021.
\newblock \url{https://www.govinfo.gov/content/pkg/DCPD-202100095/pdf/DCPD-202100095.pdf}.

\bibitem[{Fahrig(2020)}]{fahrig2020several}
Fahrig, L. 2020.
\newblock Why do several small patches hold more species than few large patches?
\newblock \emph{Global Ecology and Biogeography}, 29(4): 615--628.

\bibitem[{Faust et~al.(2016)Faust, Simmonis, Waddell, and Long}]{faust2016red}
Faust, L.~J.; Simmonis, J.~S.; Waddell, W.; and Long, S. 2016.
\newblock Red Wolf (Canis rufus) Population Viability Analysis – Report to U.S. Fish and Wildlife Service.
\newblock Technical report, Lincoln Park Zoo, Chicago.

\bibitem[{Ferraz et~al.(2021)Ferraz, Morato, Bovo, da~Costa, Ribeiro, de~Paula, Desbiez, Angelieri, and Traylor-Holzer}]{ferraz2021bridging}
Ferraz, K. M. P. M. d.~B.; Morato, R.~G.; Bovo, A. A.~A.; da~Costa, C. O.~R.; Ribeiro, Y. G.~G.; de~Paula, R.~C.; Desbiez, A. L.~J.; Angelieri, C. S.~C.; and Traylor-Holzer, K. 2021.
\newblock Bridging the gap between researchers, conservation planners, and decision makers to improve species conservation decision-making.
\newblock \emph{Conservation Science and Practice}, 3(2): e330.

\bibitem[{Ferson et~al.(2000)Ferson, Burgman, Possingham, Ball, and Andelman}]{ferson2000mathematical}
Ferson, S.; Burgman, M.; Possingham, H.; Ball, I.; and Andelman, S. 2000.
\newblock Mathematical methods for identifying representative reserve networks.
\newblock \emph{Quantitative methods for conservation biology}, 291--306.

\bibitem[{Finnegan et~al.(2021)Finnegan, Galvez-Bravo, Silveira, T{\^o}rres, J{\'a}como, Alves, and Dalerum}]{finnegan2021reserve}
Finnegan, S.~P.; Galvez-Bravo, L.; Silveira, L.; T{\^o}rres, N.; J{\'a}como, A.~A.; Alves, G.; and Dalerum, F. 2021.
\newblock Reserve size, dispersal and population viability in wide ranging carnivores: the case of jaguars in Emas National Park, Brazil.
\newblock \emph{Animal Conservation}, 24(1): 3--14.

\bibitem[{Gupta et~al.(2019)Gupta, Dilkina, Morin, Fuller, Royle, Sutherland, and Gomes}]{gupta2019reserve}
Gupta, A.; Dilkina, B.; Morin, D.~J.; Fuller, A.~K.; Royle, J.~A.; Sutherland, C.; and Gomes, C.~P. 2019.
\newblock Reserve design to optimize functional connectivity and animal density.
\newblock \emph{Conservation Biology}, 33(5): 1023--1034.

\bibitem[{Jafari and Hearne(2013)}]{jafari2013new}
Jafari, N.; and Hearne, J. 2013.
\newblock A new method to solve the fully connected reserve network design problem.
\newblock \emph{European Journal of Operational Research}, 231(1): 202--209.

\bibitem[{Jafari et~al.(2017)Jafari, Nuse, Moore, Dilkina, and Hepinstall-Cymerman}]{jafari2017achieving}
Jafari, N.; Nuse, B.~L.; Moore, C.~T.; Dilkina, B.; and Hepinstall-Cymerman, J. 2017.
\newblock Achieving full connectivity of sites in the multiperiod reserve network design problem.
\newblock \emph{Computers \& Operations Research}, 81: 119--127.

\bibitem[{Joppa and Pfaff(2009)}]{joppa2009high}
Joppa, L.~N.; and Pfaff, A. 2009.
\newblock High and far: biases in the location of protected areas.
\newblock \emph{PloS one}, 4(12): e8273.

\bibitem[{Lacy and Breininger(2021)}]{lacy2021population}
Lacy, R.~C.; and Breininger, D.~R. 2021.
\newblock Population Viability Analysis (PVA) as a platform for predicting outcomes of management options for the Florida Scrub-Jay in Brevard County.
\newblock Technical report, {NASA}.

\bibitem[{Lawler et~al.(2020)Lawler, Rinnan, Michalak, Withey, Randels, and Possingham}]{lawler2020planning}
Lawler, J.~J.; Rinnan, D.~S.; Michalak, J.~L.; Withey, J.~C.; Randels, C.~R.; and Possingham, H.~P. 2020.
\newblock Planning for climate change through additions to a national protected area network: implications for cost and configuration.
\newblock \emph{Philosophical Transactions of the Royal Society B}, 375(1794): 20190117.

\bibitem[{Leslie et~al.(2003)Leslie, Ruckelshaus, Ball, Andelman, and Possingham}]{leslie2003using}
Leslie, H.; Ruckelshaus, M.; Ball, I.~R.; Andelman, S.; and Possingham, H.~P. 2003.
\newblock Using siting algorithms in the design of marine reserve networks.
\newblock \emph{Ecological applications}, 13(sp1): 185--198.

\bibitem[{Margules, Nicholls, and Pressey(1988)}]{margules1988selecting}
Margules, C.~R.; Nicholls, A.; and Pressey, R. 1988.
\newblock Selecting networks of reserves to maximise biological diversity.
\newblock \emph{Biological conservation}, 43(1): 63--76.

\bibitem[{Marianov, ReVelle, and Snyder(2008)}]{marianov2008selecting}
Marianov, V.; ReVelle, C.; and Snyder, S. 2008.
\newblock Selecting compact habitat reserves for species with differential habitat size needs.
\newblock \emph{Computers \& Operations Research}, 35(2): 475--487.

\bibitem[{Maxwell et~al.(2020)Maxwell, Cazalis, Dudley, Hoffmann, Rodrigues, Stolton, Visconti, Woodley, Kingston, Lewis et~al.}]{maxwell2020area}
Maxwell, S.~L.; Cazalis, V.; Dudley, N.; Hoffmann, M.; Rodrigues, A.~S.; Stolton, S.; Visconti, P.; Woodley, S.; Kingston, N.; Lewis, E.; et~al. 2020.
\newblock Area-based conservation in the twenty-first century.
\newblock \emph{Nature}, 586(7828): 217--227.

\bibitem[{Morrison, Wardle, and Castley(2016)}]{morrison2016repeatability}
Morrison, C.; Wardle, C.; and Castley, J.~G. 2016.
\newblock Repeatability and reproducibility of population viability analysis (PVA) and the implications for threatened species management.
\newblock \emph{Frontiers in Ecology and Evolution}, 4: 98.

\bibitem[{{\"O}nal et~al.(2016){\"O}nal, Wang, Dissanayake, and Westervelt}]{onal2016optimal}
{\"O}nal, H.; Wang, Y.; Dissanayake, S.~T.; and Westervelt, J.~D. 2016.
\newblock Optimal design of compact and functionally contiguous conservation management areas.
\newblock \emph{European Journal of Operational Research}, 251(3): 957--968.

\bibitem[{Possingham et~al.(1993)Possingham, Day, Goldfinch, and Salzborn}]{possingham1993mathematics}
Possingham, H.; Day, J.; Goldfinch, M.; and Salzborn, F. 1993.
\newblock The mathematics of designing a network of protected areas for conservation.
\newblock In \emph{Decision Sciences: Tools for Today. Proceedings of 12th National ASOR Conference}, 536--545. ASOR Adelaide.

\bibitem[{ReVelle, Williams, and Boland(2002)}]{revelle2002counterpart}
ReVelle, C.~S.; Williams, J.~C.; and Boland, J.~J. 2002.
\newblock Counterpart models in facility location science and reserve selection science.
\newblock \emph{Environmental Modeling \& Assessment}, 7: 71--80.

\bibitem[{Rodrigues and Cazalis(2020)}]{rodrigues2020multifaceted}
Rodrigues, A.~S.; and Cazalis, V. 2020.
\newblock The multifaceted challenge of evaluating protected area effectiveness.
\newblock \emph{Nature Communications}, 11(1): 5147.

\bibitem[{Rodrigues, Orestes~Cerdeira, and Gaston(2000)}]{rodrigues2000flexibility}
Rodrigues, A.~S.; Orestes~Cerdeira, J.; and Gaston, K.~J. 2000.
\newblock Flexibility, efficiency, and accountability: adapting reserve selection algorithms to more complex conservation problems.
\newblock \emph{Ecography}, 23(5): 565--574.

\bibitem[{Rondinini and Chiozza(2010)}]{rondinini2010quantitative}
Rondinini, C.; and Chiozza, F. 2010.
\newblock Quantitative methods for defining percentage area targets for habitat types in conservation planning.
\newblock \emph{Biological Conservation}, 143(7): 1646--1653.

\bibitem[{Rosing, ReVelle, and Williams(2002)}]{rosing2002maximizing}
Rosing, K.; ReVelle, C.; and Williams, J. 2002.
\newblock Maximizing species representation under limited resources: a new and efficient heuristic.
\newblock \emph{Environmental Modeling \& Assessment}, 7: 91--98.

\bibitem[{Runge et~al.(2017)Runge, Sanders-Reed, Langtimm, Hostetler, Martin, Deutsch, Ward-Geiger, and Mahon}]{runge2017status}
Runge, M.~C.; Sanders-Reed, C.~A.; Langtimm, C.~A.; Hostetler, J.~A.; Martin, J.; Deutsch, C.~J.; Ward-Geiger, L.~I.; and Mahon, G.~L. 2017.
\newblock Status and threats analysis for the Florida manatee (Trichechus manatus latirostris), 2016.
\newblock Technical report, US Geological Survey.

\bibitem[{Schl{\"u}ter, Egea, and Banga(2009)}]{schluter2009extended}
Schl{\"u}ter, M.; Egea, J.~A.; and Banga, J.~R. 2009.
\newblock Extended ant colony optimization for non-convex mixed integer nonlinear programming.
\newblock \emph{Computers \& Operations Research}, 36(7): 2217--2229.

\bibitem[{Sinclair et~al.(2018)Sinclair, Milner-Gulland, Smith, McIntosh, Possingham, Vercammen, and Knight}]{sinclair2018use}
Sinclair, S.~P.; Milner-Gulland, E.; Smith, R.~J.; McIntosh, E.~J.; Possingham, H.~P.; Vercammen, A.; and Knight, A.~T. 2018.
\newblock The use, and usefulness, of spatial conservation prioritizations.
\newblock \emph{Conservation Letters}, 11(6): e12459.

\bibitem[{Stevens(2014)}]{stevens2014indigenous}
Stevens, S. 2014.
\newblock \emph{Indigenous peoples, national parks, and protected areas: a new paradigm linking conservation, culture, and rights}.
\newblock University of Arizona Press.

\bibitem[{Stralberg et~al.(2009)Stralberg, Applegate, Phillips, Herzog, Nur, and Warnock}]{stralberg2009optimizing}
Stralberg, D.; Applegate, D.~L.; Phillips, S.~J.; Herzog, M.~P.; Nur, N.; and Warnock, N. 2009.
\newblock Optimizing wetland restoration and management for avian communities using a mixed integer programming approach.
\newblock \emph{Biological conservation}, 142(1): 94--109.

\bibitem[{Underhill(1994)}]{underhill1994optimal}
Underhill, L. 1994.
\newblock Optimal and suboptimal reserve selection algorithms.
\newblock \emph{Biological Conservation}, 70(1): 85--87.

\bibitem[{Watson et~al.(2014)Watson, Dudley, Segan, and Hockings}]{watson2014performance}
Watson, J.~E.; Dudley, N.; Segan, D.~B.; and Hockings, M. 2014.
\newblock The performance and potential of protected areas.
\newblock \emph{Nature}, 515(7525): 67--73.

\bibitem[{West, Igoe, and Brockington(2006)}]{west2006parks}
West, P.; Igoe, J.; and Brockington, D. 2006.
\newblock Parks and peoples: the social impact of protected areas.
\newblock \emph{Annu. Rev. Anthropol.}, 35: 251--277.

\bibitem[{Williams, Rondinini, and Tilman(2022)}]{williams2022global}
Williams, D.~R.; Rondinini, C.; and Tilman, D. 2022.
\newblock Global protected areas seem insufficient to safeguard half of the world's mammals from human-induced extinction.
\newblock \emph{Proceedings of the National Academy of Sciences}, 119(24): e2200118119.

\bibitem[{Williams et~al.(2020)Williams, Scriven, Burslem, Hill, Reynolds, Agama, Kugan, Maycock, Khoo, Hastie et~al.}]{williams2020incorporating}
Williams, S.~H.; Scriven, S.~A.; Burslem, D.~F.; Hill, J.~K.; Reynolds, G.; Agama, A.~L.; Kugan, F.; Maycock, C.~R.; Khoo, E.; Hastie, A.~Y.; et~al. 2020.
\newblock Incorporating connectivity into conservation planning for the optimal representation of multiple species and ecosystem services.
\newblock \emph{Conservation Biology}, 34(4): 934--942.

\bibitem[{Winton, Bishop, and Larsen(2020)}]{winton2020protected}
Winton, S.~A.; Bishop, C.~A.; and Larsen, K.~W. 2020.
\newblock When protected areas are not enough: low-traffic roads projected to cause a decline in a northern viper population.
\newblock \emph{Endangered Species Research}, 41: 131--139.

\bibitem[{{World Wildlife Fund}(1980)}]{world1980world}
{World Wildlife Fund}. 1980.
\newblock \emph{World conservation strategy: Living resource conservation for sustainable development}, volume~1.
\newblock Gland, Switzerland: {IUCN}.

\end{thebibliography}
